\documentclass[a4paper,10pt]{article}
\usepackage[utf8x]{inputenc}

\usepackage{ulem}
\usepackage{cancel}

\usepackage{minitoc}
\newcommand{\nomi}{\mathbf{i}}
\newcommand{\nomj}{\mathbf{j}}
\newcommand{\nomk}{\mathbf{k}}

\newcommand{\unomj}{\underline{\mathbf{j}}}

\newcommand{\cnomm}{\mathbf{m}}
\newcommand{\cnomn}{\mathbf{n}}

\newcommand{\ucnomm}{\underline{\mathbf{m}}}

\usepackage{multicol}
\usepackage{bussproofs}
\usepackage{wrapfig}

\usepackage[english]{mytheorems} 
\usepackage{amsmath,amsthm,amssymb}
\usepackage{stmaryrd} 

\usepackage{enumitem}

\usepackage{graphicx}
\usepackage[all]{xy}
\usepackage{hyperref}
\usepackage{color}
\usepackage{float}
\graphicspath{{Images/}}

\usepackage{fullpage}

\usepackage{latexsym}
\usepackage{txfonts}

\usepackage{color}

\newcommand{\D}{D^+}
\newcommand{\set}[1]{\{\,#1\,\}}

\newcommand{\M}{\mathcal{M}}

\RequirePackage{graphicx}

\RequirePackage{xspace}

\newcommand{\mfLL}{\mf{L}_L}

\newcommand{\X}{\mathcal{X}}
\renewcommand{\P}{\mathcal{P}} %
\newcommand{\rto}{\xrightarrow{\;\;\;}}
\newcommand{\rTo}{\rto}

\newcommand{\tiff}{\text{ iff }}
\newcommand{\tif}{\text{ if }}
\newcommand{\tand}{\text{ and }}

{ \begin{list}%
	{$\bullet$}%
	{\setlength{\labelwidth}{30pt}%
	 \setlength{\leftmargin}{35pt}%
	 \setlength{\itemsep}{\parsep}}}%
{ \end{list} }

{\begin{list}%
    {{$(\roman{enumi})$}\stepcounter{enumi}}%
    {\setcounter{enumi}{1}%
      \setlength{\labelwidth}{0pt}%
      \setlength{\leftmargin}{5pt}%
      \setlength{\listparindent}{10pt}%
      \setlength{\labelsep}{5pt}%
    }}%
{\end{list}}

\newcommand{\mc}{\mathcal}
\newcommand{\mf}{\mathfrak}
\newcommand{\mb}{\mathbb}
\newcommand{\mbf}{\mathbf}
\newcommand{\msf}{\mathsf}
\newcommand{\llb}{\llbracket}
\newcommand{\rrb}{\rrbracket}
\newcommand{\lb}{\lbrack}
\newcommand{\rb}{\rbrack}

\newcommand{\W}{\wedge}
\newcommand{\V}{\vee}
\newcommand{\ef}{ ( \exists \forall ) }
\newcommand{\fe}{ ( \forall \exists ) }

\newcommand{\cov}{\lhd}

\newcommand{\covD}{\langle \cov \rangle}
\newcommand{\covB}{\lb \cov \rb}

\newcommand{\inD}{\langle \in \rangle}

\newcommand{\niD}{\langle \ni \rangle}
\newcommand{\niB}{\lb \ni \rb}

\newcommand{\tor}{\text{ or }}

\newcommand{\ackRcl}{RA\cl}

\newcommand{\MT}{MT}
\newcommand{\DeMorgan}{DM}
\newcommand{\tra}{TR}

\newcommand{\F}{\mb{F}}

\newcommand{\aj}{AJ}
\newcommand{\ajW}{RS\wedge}
\newcommand{\ajV}{RS\vee}

\newcommand{\ajB}{AJ\square}

\newcommand{\MinCovB}{MinCov2}

\newcommand{\ap}{AP}
\newcommand{\apD}{AP\diamondsuit}

\newcommand{\tW}{T\wedge}
\newcommand{\tV}{T\vee}
\newcommand{\tminus}{T\smallsetminus}
\newcommand{\tWB}{T\W\bot}
\newcommand{\tNN}{TNN}

\newcommand{\tnm}{TNM}

\newcommand{\baW}{BA\wedge}
\newcommand{\baV}{BA\vee}

\newcommand{\tBD}{TBD}
\newcommand{\tDB}{TDB}
\newcommand{\tRRminus}{TRR^{-1}}

\newcommand{\spW}{SP\wedge}
\newcommand{\spV}{SP\vee}

\newcommand{\AtCoat}{AtCoat1}
\newcommand{\AtCoattwo}{AtCoat2}
\newcommand{\Vbot}{\vee\bot}
\newcommand{\Wtop}{\wedge\top}

\newcommand{\leqAtom}{Atom\rxx}

\newcommand{\distWV}{D\wedge\vee}
\newcommand{\distVW}{D\vee\wedge}
\newcommand{\asV}{A\vee}
\newcommand{\asW}{A\wedge}
\newcommand{\cmV}{C\vee}
\newcommand{\cmW}{C\wedge}

\newcommand{\bis}{bis}

\newcommand{\Substitution}{Sub}

\RequirePackage{amsmath,amssymb}
\RequirePackage{eufrak}

\newcommand{\tobesaid}[1]{%
  \ifx\undefined\public
  \begin{itemize}
  \item[\texttt{>>}] \texttt{#1} 
  \end{itemize}
  \else\fi
}

\newcommand{\margtobesaid}[1]{%
  \ifx\undefined\public
  \marginpar{\texttt{#1}}%
  \else\fi
}

\newcommand{\Powerset}{\mathcal{P}}
\renewcommand{\P}{\Powerset}

\newcommand{\ud}[1]{\underline{#1}}

\newcommand{\rxy}{R_{XY}}
\newcommand{\ryx}{R_{YX}}
\newcommand{\rxx}{R_{XX}}

\newcommand{\Drxy}{\langle \rxy \rangle}
\newcommand{\Brxy}{\lb \rxy \rb}
\newcommand{\Drxym}{\langle \rxy^{-1} \rangle}
\newcommand{\Brxym}{\lb \rxy^{-1} \rb}

\newcommand{\Dryx}{\langle \ryx \rangle}
\newcommand{\Bryx}{\lb \ryx \rb}
\newcommand{\Dryxm}{\langle \ryx^{-1} \rangle}
\newcommand{\Bryxm}{\lb \ryx^{-1} \rb}

\newcommand{\Drxx}{\langle \rxx \rangle}

\newcommand{\leqx}{\leq}
\newcommand{\leqy}{\preceq}

\newcommand{\leqjD}{\langle \leq_J \rangle}
\newcommand{\leqxD}{\langle \leq_X \rangle}

\newcommand{\sys}[1]{\left%
    (\begin{array}[c]{r@{\hspace{2pt}}l}#1\end{array}\right )}

\renewcommand{\implies}{\Rightarrow}

\newcommand{\false}{\mathtt{false}}

\newcommand{\fD}{\mc{D}}

\newcommand{\A}{\mb{A}}
\renewcommand{\X}{\mb{X}}

\newcommand{\cl}{cl}
\newcommand{\clPres}{\overline{cl}}

\newcommand{\AtProp}{\msf{AtProp}}
\newcommand{\Var}{\msf{Var}}
\newcommand{\Nom}{\msf{Nom}}
\newcommand{\CNom}{\msf{CNom}}

\newcommand{\jty}{J^{\infty}}
\newcommand{\mty}{M^{\infty}}

\allowdisplaybreaks[1]

\title{Dual characterizations for finite lattices via correspondence theory for monotone modal logic}
\author{Sabine Frittella, Alessandra Palmigiano and Luigi Santocanale}
\date{}
\begin{document}

\maketitle

\begin{abstract}
    We establish a formal connection between algorithmic correspondence theory and certain dual characterization 
    results for finite lattices,
    similar to 
    Nation's
characterization of a hierarchy of pseudovarieties of finite lattices, progressively generalizing finite distributive lattices. This formal connection is mediated through monotone modal logic.
Indeed, we adapt the correspondence algorithm ALBA to the setting of monotone modal logic, and we use a certain duality-induced encoding of finite lattices as monotone neighbourhood frames to translate lattice terms into formulas in monotone modal logic. 


\end{abstract}

\tableofcontents

\section{Introduction}

\paragraph*{Dual characterization results for finite lattices.}
The present paper builds on a duality for finite lattices, established
by  Santocanale  \cite{duality}. The structures dually
equivalent to finite lattices are referred to
as {\it join-presentations}, and are certain triples $(
X,\leq,\mc{M})$ such that $(X,\leq)$ is a finite
poset, and $\mc{M} : X \rTo \mathcal{P}\mathcal{P}X$. In \cite{duality}, it has been pointed
out---and indicated as a worthwhile research direction---that the
existence of this duality makes it possible to investigate systematic
{\it dual characterization} results, between equations or inequalities in the algebraic language of
lattices on one side, and first-order conditions in the language
of join-presentations on the other.
One
significant instance of such systematic dual characterizations has
been developed in the same paper, between a class of inequalities in
the language of lattices and a corresponding class of first-order conditions. Both classes are  parametric in the class of finite trees (cf.\
\cite[Proposition 8.5]{duality}). This result generalizes Nation's
\cite[Section 5]{NationApproach}
 stating that a certain  class of finite lattices\footnote{Namely, the finite lattices such that the
length of their $D$-{\it chains} has a uniform upper bound.}
is a pseudovariety, and is similar to
Semenova's results \cite{SemenovaSuborders}.


\paragraph*{From modal logic to unified correspondence theory.}
Modal logic is an area in which systematic dual characterization results have been extensively developed, giving rise to a very rich theory---the so called {\it modal correspondence theory}---which has been investigated for almost forty years. Modal correspondence theory was originally developed in a purely model-theoretic way \cite{vB83}.
However,  correspondence-related phenomena have been studied in an algebraic framework subsuming duality theory since the early 90s \cite{JonssonTarski51}, and very recently, a unified correspondence framework has emerged \cite{unifiedCorrespondence}, which is based on duality, and uniformly extends correspondence theory to many nonclassical logics.
One of the main tools developed by this theory is an algorithm (actually various cognate versions of it, cf.\ \cite{ALBA}) or calculus for correspondence, called ALBA, which mechanizes dual characterization meta-arguments. In particular, as discussed in \cite{ALBA} and \cite{Gilardi}, the core of ALBA is the encoding of a general meta-argument for correspondence known in the literature as the \textit{minimal valuation argument} into a rule which relies on the\textit{ Ackermann lemma} \cite{Ackermann:Untersuchung}. The algorithm ALBA takes in input formulas or inequalities in a given propositional language and, whenever it succeeds\footnote{It is well known \cite{Chagrov:Chagrova:2006} that the problem of whether a formula admits a first-order correspondent is undecidable.}, it computes the {\it first-order correspondent} of the given formula or inequality, i.e.\ a first-order sentence which holds in a given structure exactly when the given propositional formula or inequality is valid in the dual algebra of that structure. The general theory also provides the syntactic characterization of a class of formulas/inequalities for each logic, the so-called inductive formulas/inequalities, on which the algorithm is guaranteed to uniformly succeed.  For each language, inductive inequalities form the largest  such class syntactically defined so far in the literature.

\paragraph*{Aim of the paper.}
Given the availability of this theory, it seems natural to try and understand dual characterization results such as \cite[Proposition 8.5]{duality} as instances of a more general unified correspondence mechanism. 
This is what the present paper aims at doing,
by establishing a novel dual characterization result similar to Nation's. Our result paves the way to the mechanization and systematization of dual characterizations  such as the one in \cite{duality}.

\paragraph*{Methodology: basic algorithmic correspondence for monotone modal logic.}
Our approach is based on an adaptation of the algorithm/calculus ALBA of \cite{ALBA} to the case of {\it monotone modal logic}. This adaptation is necessary, since some of the rules in the standard version of the algorithm would not be sound for the modal connectives of monotone modal logic, and is one of the contributions of the present paper.
The adapted ALBA is semantically justified in the general environment of two-sorted frames (cf.\ Section \ref{sec : 2sorted}),
which are general structures that can encode monotone neighbourhood frames as special cases.
As their name suggests, two-sorted frames are relational structures
based on two domains. Normal modal operators can be associated
in the standard way
with the binary relations on two-sorted frames.
Monotone modal operators can be then interpreted on two-sorted frames as the composition of some of these normal modalities.
This provides the basic semantic environment for the adapted ALBA.

Correspondence theory for monotone modal logic has already been studied in \cite{HelleThesis}, where a class of monotone modal formulas which are guaranteed to have a first order correspondent has been identified.
However, the class of inductive inequalities corresponding to the ALBA setting is strictly larger than the one in \cite{HelleThesis}.

\paragraph*{Enhancing ALBA.}
However, the translations of inequalities such as Nation's \cite{NationApproach} and as the ones treated in the present paper
fall outside the inductive class. Hence, another contribution of the present paper is the addition of special rules which are sound on the specific semantic setting arising from finite lattices.
Interestingly, an Ackermann-type rule features among these additional rules, the soundness of which cannot be straightforwardly explained in terms of the Ackermann lemma, but which however still intuitively encodes a minimal valuation argument.


%
%


\paragraph*{Organization of the paper.}
In Section 2, we collect  preliminaries on
the duality between finite lattices and join-presentations,
the language and neighbourhood semantics of monotone modal logic,
the duality-induced `standard translation' of lattice terms into terms in the language of monotone modal logic,
and the algorithm for correspondence ALBA.
In Section 3, we adapt the algorithm ALBA  specifically to monotone modal logic via the introduction of two-sorted frames.
In Section 4, we enhance the adapted ALBA by introducing additional rules, and prove their soundness w.r.t.\
the semantic environment of
so-called enriched two-sorted frames which can be naturally associated with finite lattices.
In Section 5,  upper bounds on the length of $\D$-chains (cf.\ Definition \ref{def : Dchain}) are  obtained as a reduction of the enhanced ALBA. 
Section 6 collects the conclusions and further directions.
The proof of a technical lemma appears in the appendix.

\section{Preliminaries}
\label{sec : Preliminaries}
    The aim of the present section is to collect preliminaries
belonging to diverse fields of logic, and to connect them
so as to set the stage for the main result.
In the next subsection, we report on a duality on objects for finite lattices, which has been 
introduced in \cite{duality}. 
Given that the structures dual to finite lattices can be naturally associated with monotone neighbourhood frames, and given that monotone neighbourhood frames are standard models for monotone modal logic, 
the duality presented in subsection \ref{subsec:join-presentation}
serves as a basis for the definition of a standard translation 
between lattice terms and monotone monotone modal logic formulas.
In subsection \ref{subsec : MML}, we recall the basic definitions about monotone modal logic and neighbourhood frames, and we show how we can represent a finite lattice as a monotone neighbourhood frame.
In subsection \ref{subsec : standard translation}, we define a standard translation between lattice terms and formulas of monotone modal logic, and show that this translation adequately transfers and reflects the validity of lattice inequalities on any finite lattice $L$, and  the validity of their  standard translations on the monotone neighbourhood frames associated with $L$. 
Finally, in subsection \ref{subsec : intro ALBA}, we give an informal    presentation of the algorithm for correspondence ALBA, introduced in \cite{ALBA}, for correspondence for normal modal logic.

\subsection{Dual equivalence for finite lattices}
\label{subsec:join-presentation}
In the present subsection, we report on the object-part of a dual
equivalence between finite lattices and certain poset-based structures
(cf.\ Definition \ref{def : join presentation}). Our presentation is
based on \cite{duality}.  These structures will turn out to be special
neighbourhood frames, and hence the existence of this duality provides
the bridge between the propositional logic of lattices and monotone
modal logic.

In what follows, $L$ will denote a finite lattice. Elements of $L$ will be denoted  $a, b \ldots$ Throughout the paper, the letters  $i, j, k$ will be reserved for join-irreducible elements of $L$ (the set of which is $J(L)$), and $m, n$ for meet-irreducible elements of $L$ (the set of which is $M(L)$), respectively.
Recall that an element $j\neq \bot$ of $L$ is \textit{join-irreducible} iff   $j = a \vee b$ implies that either $j=a$ or $j=b$ for all $a, b\in L$. Order-dually, an element $m\neq \top$ of $L$ is  \textit{meet-irreducible} iff $m = a \wedge b$ implies that either $m=a$ or $m=b$ for all $a, b\in L$.
A subset $C\subseteq L$ is a \textit{join-cover} of $a\in L$ if $a \leq \bigvee C$.

For any poset $(S,\leq)$, its associated \textit{refinement relation}, denoted $\ll$, is defined on the set $\mc{P}_{f}(S)$ of finite subsets of $S$ by the following stipulation:

\begin{equation}
\label{eq : def refinement}
A \ll B\ \quad \mbox{ iff }\ \quad \text{ for every } a \in A \; \text{ there exists some } b \in B \mbox{ such that } a \leq b.
\end{equation}
Equivalently, 
\[
A \ll B\ \quad \mbox{ iff }\ \quad \mbox{ iff }\ \quad {\downarrow} A \subseteq\: {\downarrow} B,
\]
where 
${\downarrow} C := \{ x \in S \mid x\leq c  $ for some $c$ in $C \}$ for every $C \subseteq S$.
Throughout the paper, 
we say that a join-cover $C$ of $a$ is \textit{minimal} 
if it is an $\leq$-antichain, and if, for any $\leq$-antichain $D\subseteq L$, ($a\leq \bigvee D$ and  $D \ll C$) imply $D=C$.
A join-cover of $a$ is \textit{trivial} if it contains $a$. We can easily show that  any join-irreducible element $j\in J(L)$ has only one trivial minimal join-cover, that is, the singleton $\{ j\}$.

\paragraph*{Direct presentations and their closure operators.}
In the present    paragraph, we define direct presentations, and
introduce a closure operator over these presentations which 
is a key ingredient of the duality on objects between finite lattices and reflexive and transitive presentations (see paragraph below).

\begin{definition}
    A \textit{presentation} is a triple $( X, \leq, \mc{M} )$ such that $(X, \leq)$ is a poset, and $\mc{M}: X \to \mc{P}\mc{P}X$. A presentation is 
    		\begin{itemize}
			\item \textit{monotone} if for all $ x,y \in X$,  any $ C \subseteq  X$,
			if $ y \leq x$ and $C \in \M(x)$,
			then  $D\ll C$ for some $ D \in \M(y)$;
			\item \textit{reflexive} if for each $ x \in X$, there exists some  $ C \in \M(x)$ such that $C \ll \{x\}$;
			\item \textit{transitive} if for every $x\in X$ and every $C \subseteq X$, 
			if $C \in \M(x)$ 
			then for
			every collection $\{ D_c \mid c\in C\}$ such that $D_c \in \M(c)$ for every $c \in C$,
			there exists some $E\in \M(x)$
			such that  $E \ll \bigcup_{c\in C} D_{c}$;
			\item \textit{direct} if it is monotone, reflexive, and transitive.
		\end{itemize}
\end{definition}

Recall that a \textit{downset} of $(X,\leq)$ is a subset $S\subseteq X$
 such that for all $ x,y \in X$, if $ y\leq x$ and $x\in S$ then $y\in S$. Let $\fD (X,\leq)$ denote the set of downsets of $(X,\leq)$.

For any  presentation $\A = (\X, \M)$ where $\X:=(X,\leq)$ is a poset, the assignment $\clPres_\A : \fD \X \longrightarrow \fD \X$ is defined  as follows: 
for any $S \in \fD \X $,  
\begin{equation}
\label{eq : def closure}
\clPres_\A(S) := \{ x\in X \mid D\subseteq S \text{ for some } D\in \M(x)\} .
\end{equation}

\begin{lemma}
For any direct presentation $\A = (X,\leq,\M)$, the map $\clPres_\A$
is a closure operator.
\end{lemma}
\begin{proof}
We first prove that the map $\clPres_\A$ is well-defined. 
Fix $S\in \fD \X$, $x \in \clPres_\A(S)$ and $y\in X$. Assume that $y\leq x$. 
Since $x\in \clPres_\A(S)$, there is some $C_x \in \M(x)$ such that $C_x \subseteq S$.
In addition, 
since $\A$ is monotone, 
$y\leq x$ implies that 
there exists some $C_y \in \M(y)$ such that 
$C_y \ll C_x $. By definition of $\ll$, we have that $C_y\subseteq S$
because $S$ is a downset.
 Thus $y\in \clPres_\A(S)$. This finishes the proof that $\clPres_\A(S)$ is a downset. 
Hence the map $\clPres_\A$ is well-defined.

To prove that $\clPres_\A$ is  a closure operator, we need to show that 
$\clPres_\A$ is order-preserving, and that
$S \subseteq \clPres_\A (S) $ and $\clPres_\A (\clPres_\A (S))\subseteq \clPres_\A (S)$
for any $S\in \fD\X$. It is immediate to see that $\clPres_\A$ is order preserving. 

Since $\A$ is reflexive,  there is some $C\in \M(x)$ such that  $C \subseteq {\downarrow}x$. Moreover, $x\in S$ implies that  ${\downarrow}x \subseteq S$.
Hence, by definition of $\clPres_\A$, we have that $x\in \clPres_\A(S)$ for any $x\in S$, that is $S\subseteq \clPres_\A(S)$.

It remains to be shown that $\clPres_\A (\clPres_\A (S))\subseteq \clPres_\A (S)$
for any $S\in \fD\X$. Let  $x \in \clPres_\A (\clPres_\A (S))$. 
By definition of $\clPres_\A$, there exists some $D \in \M(x)$ such that $D\subseteq \clPres_\A(S) $. 
Then any $d\in D$ is an element of 
$\clPres_\A(S)$. Thus, 
for each $d\in D$  there exists some $E_d \in \M(d)$ such that $E_d\subseteq S $. 
Since $\A$ is transitive, there is some $C\in \M(x)$, such that   $C \ll \bigcup_{d\in D} E_d $. Thus $C \subseteq S$, and, by definition of $\clPres_\A(S)$, this proves that $x\in \clPres_\A(S)$. 
This completes the proof that $\clPres_\A$ is a closure operator.
\end{proof}

\begin{definition}
    For any direct presentation $\A = ( X,\leq ,\M )$,
    a downset $S \subseteq X$ is \textit{closed} if
    $S = \clPres_\A(S)$. The closure of a downset $S\subseteq X$ is the set $\clPres_\A(S)$. In the following, whenever it causes no confusion, we denote the closure of a downset $S$ by $\overline{S}$. 
\end{definition}

Notice that for any direct presentation $\A = ( X,\leq ,\M )$, we can  extend the closure operator $\clPres_\A$ to sets, as follows:
\begin{align}
\label{eq : def closure on sets}
\cl_\A : \P X & \longrightarrow \P X  \notag\\
S & \longmapsto \clPres_\A({\downarrow_{\leq}} S).
\end{align}
Since, $\clPres_\A$ and $ {\downarrow_{\leq}}$ are closure operators on downsets and on sets respectively, we can easily prove the $\cl_\A$ is a closure operator too.

\paragraph*{Join-presentation of a finite lattice.}

\begin{definition}
\label{def : join presentation}
	The \textit{join-presentation}\footnote{Join-presentations are also referred to as $OD$-graphs in the literature (cf.\ \cite{NationApproach, duality}).}
	of a lattice $L$ is the presentation $( J(L), \leq, \mc{M} )$ such that
	 $(J(L),\leq)$ is the  poset of the join-irreducible elements of $L$ with the order induced by $L$,
	 and $\mc{M}$ is the map $J(L) \longrightarrow \mc{P}\mc{P}J(L)$
	assigning any $j$ to the collection of its minimal join-covers.
\end{definition}

\begin{lemma}[cf.\ Lemma 4.2 in \cite{duality}]
For any finite lattice $L$, 
the join-presentation $(J(L),\leq, \M) $ associated with $L$  is a direct presentation.
\end{lemma}

More generally, we can associate every element $a$ of a lattice $L$
with the set $\mc{M}(a)$ of its minimal join-covers.
The following lemma lists some properties
of $\mc{M} : L \longrightarrow \P\P J(L)$.
\begin{lemma}[cf.\ \cite{duality}, page 5]
\label{lem : minimal covers}
    Let $(L, \leq)$ be a finite lattice. For all $a \in  L$, $j\in J(L)$, $C\in \M(a)$ and $Y\subseteq L$,
    \begin{enumerate}
    \item $C \subseteq J(L)$, and $C$ is an $\leq$-antichain;
    \item $\mc{M}(a)$ is a $\ll$-antichain;
     \item if $a \leq \bigvee Y$,
     then there exists some $D \in \mc{M}(a)$ such that
     $D \ll Y$;
     \item $\{ \,j\, \} \in \mc{M}(j)$.
    \end{enumerate}
\end{lemma}

    For every finite lattice $L$,
let $\mfLL$ be the lattice of the closed downsets
    of the join-presentation $(J(L),\leq, \M )$ associated with $L$.

\begin{proposition}[cf.\ \cite{NationApproach}]
 Every finite lattice $L$ is
isomorphic to the lattice $\mfLL$ as above.
\end{proposition}

The following lemmas will be useful in the remainder of the paper.

\begin{lemma}[Lemma 4.2 in \cite{duality}]
\label{lem : closure down j}
For any finite lattice $L$ and any $j\in J(L)$, 
the downset  ${\downarrow_{J(L)}}j$ is a closed
subset of the join-presentation $(J(L),\leq, \M)$ associated with $L$.
\end{lemma}

\begin{lemma}
\label{lem : minimal cover closure}
Let $L$ be a finite lattice, and  $(J(L),\leq, \M)$ be  its join-presentation. For any $j,k \in J(L)$ and any $\leq$-antichain $C \subseteq J(L)$, if $C \in \M(j)$ and $k \in C$,
\begin{enumerate}
\item$j \notin \overline{{\downarrow_{\leq}}( C \smallsetminus k)}$,
\item$k \notin \overline{{\downarrow_{\leq}}( C \smallsetminus k)}$,
\item $j \notin \{ k'\in J(L) \mid k' < k\}$, 
\item there is no $D \in \M(j)$ such that $D \subseteq \overline{{\downarrow_{\leq_J}}(C \smallsetminus k)} \cup \{ k'\in J(L) \mid k' < k\}$.
\end{enumerate}
\end{lemma}
\begin{proof}Fix $C \in \M(j)$ and $k\in C$.
As to item 1.  Since $C$ is a minimal cover of $j$, the sets $C\smallsetminus k$ and ${\downarrow_{\leq_J}}(C\smallsetminus k)$ are not  covers of $j$. Hence $j \notin \overline{{\downarrow_{\leq}}( C \smallsetminus k)}$.

We show item 2 by contradiction. Assume that $k \in \overline{{\downarrow_{\leq}}( C \smallsetminus k)}$. By the definition of  closure, this implies that
there exists some $D \in \M(k)$ such that $D \subseteq {\downarrow_{\leq}}( C \smallsetminus k)$. 
The following chain of inequalities holds
\begin{align*}
& \; j \leq \bigvee  C  
\tag{$C\in \M(j)$}\\
=& \;  \bigvee ((C\smallsetminus k) \cup \{ k\})
\\
= & \; (\bigvee (C\smallsetminus k) )\vee  k \\
\leq & \; \bigvee ( C \smallsetminus k) \vee  \bigvee D
\tag{$D\in\M(k)$} \\
= & \; \bigvee (( C\smallsetminus k) \cup D ),
\end{align*}
which shows that
the set $( C \smallsetminus k) \cup D$ is a cover of $j$. Hence, there exists a minimal cover $C'\in \M(j)$ that refines it, i.e.\ such that  $C' \ll ( C\smallsetminus k ) \cup D$. 
By the definition of $\ll$, this means that $C' \subseteq {\downarrow_{\leq}} (( C\smallsetminus k ) \cup D) $. Since 
$D \subseteq {\downarrow_{\leq}}( C \smallsetminus k)$, we have that ${\downarrow_{\leq}} (( C\smallsetminus k ) \cup D) =  {\downarrow_{\leq}}( C \smallsetminus k)$, which proves  that $C' \subseteq {\downarrow_{\leq}}( C \smallsetminus k)$.  This proves that $j \in \overline{{\downarrow_{\leq}}( C \smallsetminus k)}$, which contradicts item 1.

Item 3 immediately follows from the definition of a minimal cover.

As to item 4, suppose for contradiction that there exists some  $D \in \M(j)$ such that 
$D \subseteq \overline{{\downarrow_{\leq_J}}(C \smallsetminus k)} \cup \{ k'\in J(L) \mid k' < k\}$.
Then, for any $d\in D$, there exists some $k_d \in \overline{{\downarrow_{\leq_J}}(C \smallsetminus k)} \cup \{ k'\in J(L) \mid k' < k\}$ such that  $d\leq k_d$. If $k_d \in \{ k'\in J(L) \mid k' < k\}$, then $k_d < k$. If $k_d \notin \{ k'\in J(L) \mid k' < k\}$, then $k_d \in  \overline{{\downarrow_{\leq_J}}(C \smallsetminus k)}$ and there is some $E_d \in \M(d)$ such that 
$E_d \ll C \smallsetminus k$. Since the join-presentation of $L$ is a transitive presentation, the set 
$$ E := \bigcup \{ E_d \mid d\in D \tand k_d \notin \{ k'\in J(L) \mid k' < k\} \} \cup \bigcup  \{ k_d \mid d\in D \tand k_d < k\} $$
is a cover of $j$. 
Hence, there exists some $E' \in \M(j)$  such that $E' \ll E$. Since $E \ll C$ and the relation $\ll$ is transitive, this implies that $E'\ll C$.
Hence, to finish the proof, it is enough  to show that $E' \neq C$, which would contradict the minimality of $C$.  Since
$$
k \notin \overline{{\downarrow_{\leq_J}}(C \smallsetminus k)} \cup \{ k'\in J(L) \mid k' < k\}
\quad \tand \quad {\downarrow_{\leq_J}} E' \subseteq {\downarrow_{\leq_J}} E \subseteq \overline{{\downarrow_{\leq_J}}(C \smallsetminus k)} \cup \{ k'\in J(L) \mid k' < k\},$$ we have that  $k\notin E'$. 
Since, by assumption, $k\in C$, this proves that $E'\neq C$ as required.
\end{proof}


\subsection{An environment for correspondence}
\label{subsec : MML}
The structures described in the previous subsection
are very close to neighbourhood frames (we will expand on this at the end of the present subsection). Neighbourhood frames are well known 
to provide a state-based semantics for monotone modal logic (see \cite{HelleThesis}).
Hence, as discussed in \cite{unifiedCorrespondence},
the duality between lattices and join-presentations
induces a correspondence-type relation
between the propositional language and logic of lattices,
and a fragment of the language of monotone modal logic.

In the present section we collect the
basic ingredients of this correspondence:
the languages, their interpretations,
and a syntactic translation which may be regarded as a kind of standard translation
between the language of lattices
and the monotone modal language.

\begin{definition}
The language of lattice terms $\mc{L}_{Latt}$ over the set of
        variables $\AtProp$ is as usual given by the following syntax
        $$ \varphi ::= \bot \mid \top
        \mid p \mid \varphi \vee \varphi \mid
        \varphi \wedge \varphi, $$
        with $p \in \AtProp$.
\end{definition}

\begin{definition}
\label{def : language MML}
         The language of monotone modal logic $\mc{L}_{MML}$
        over the set of variables $\AtProp$
        is recursively defined as follows:
        $$ \varphi ::= \bot \mid \top \mid p \mid
        \neg \varphi \mid
        \varphi \vee \varphi \mid \varphi \wedge \varphi \mid
        ( \exists \forall ) \varphi \mid
        ( \forall \exists ) \varphi. $$
\end{definition}

\begin{definition}  
\label{def : neighbourhood frames}
A \textit{neighbourhood frame} is a tuple
        $\F = ( X,\sigma  )$ such that $X$ is a set and $\sigma : X \longrightarrow \P\P X$ is a map.
         For any $x\in X$, any element $N\in  \sigma(x)$ is called a \textit{neighbourhood} of $x$.
        A  neighbourhood frame $\F$ is \textit{monotone} if for any $x\in X$,
        the collection $\sigma(x)$ is an upward closed subset of $(\P X , \subseteq)$.
         A \textit{neighbourhood model} is a tuple
        $\mb{M}=( \F, v)$
        such that $\F = ( X,\sigma )$ is a neighbourhood frame
        and $v: \AtProp \longrightarrow \P X$ is a valuation.
\end{definition}

\begin{definition}
\label{def : semantic MML}
         For any neighbourhood model $\mb{M} = ( \F, v )$ and any $w\in X$, the \textit{satisfaction} of any formula  $\varphi \in \mc{L}_{MML}$ in $\mb{M}$ at $w$ is defined recursively as follows:
        \begin{align*}
        \mb{M},w \Vdash \bot &\quad\phantom{\text{ iff }} \quad never\\
        \mb{M},w \Vdash \top &\quad\phantom{\text{ iff }} \quad always \\
        \mb{M},w \Vdash p &\quad\text{ iff }\quad w \in v(p) \\
        \mb{M},w \Vdash \neg \varphi &\quad\text{ iff }\quad
        \mb{M},w \nVdash \varphi \\
        \mb{M},w \Vdash \varphi \vee \psi &\quad\text{ iff }\quad
        \mb{M},w \Vdash \varphi \text{ or }
        \mb{M},w \Vdash \psi \\
        \mb{M},w \Vdash \varphi \wedge \psi &\quad\text{ iff }\quad
        \mb{M},w \Vdash \varphi \text{ and }
        \mb{M},w \Vdash \psi \\
        \mb{M},w \Vdash ( \exists \forall ) \varphi
        &\quad\text{ iff }\quad \text{ there exists some }
         C \in \sigma(w) \text{ such that, for each } c\in C, \text{ we have }  \mb{M},c \Vdash \varphi \;
        \\
        \mb{M},w \Vdash ( \forall \exists ) \varphi
        &\quad\text{ iff }\quad
        \text{ for each } C \in \sigma(w)  \text{ there exists  some }  c\in C  \text{ such that we have }\mb{M},c \Vdash \varphi .
        \end{align*}
The above definition of local satisfaction  naturally extends 
to \textit{global satisfaction} as follows: for any formula $\varphi \in \mc{L}_{MML}$,
$$
\mb{M} \Vdash \varphi \quad \tiff \quad \mb{M},w \Vdash \varphi \text{ for any } w \in X.
$$
The notions of \textit{local} and \textit{global validity} are defined as follows:
for any formula $\varphi \in \mc{L}_{MML}$,  any neighbourhood frame $\F = (X,\sigma)$, and any $w\in X$,
$$ \F, w \Vdash \varphi \quad \tiff \quad (\F,v),w \Vdash \varphi \text{ for any valuation } v : \AtProp \longrightarrow X.$$ 
$$ \F \Vdash \varphi \quad \tiff \quad (\F,v) \Vdash \varphi \text{ for any valuation } v : \AtProp \longrightarrow X.$$ 
All the above definitions of satisfaction and validity can be naturally extended to $\mc{L}_{MML}$-inequalities as follows: for all formulas $\varphi, \psi \in \mc{L}_{MML}$, and any neighbourhood model $\mb{M} = (\F,v)$,
$$ \mb{M} \Vdash \varphi\leq\psi \quad \tiff \quad \mb{M},w \Vdash \varphi \text{ implies } \mb{M},w \Vdash \psi \text{ for any } w \in X.$$ 
$$ \F \Vdash \varphi\leq\psi \quad \tiff \quad \text{ for any valuation } v  \text{ and any } w \in X, \text{ if } (\F,v), w \Vdash \varphi \text{ then } (\F,v), w\Vdash \psi.$$ 

\end{definition}
\begin{remark}
\label{rem : semantic MML}
We notice that the definition above is usually adopted only for monotone neighbourhood frames and not for arbitrary neighbourhood frames. Under this definition, any neighbourhood frame behaves like a monotone one.
Adopting this definition, rather than the usual one, 
is more advantageous for the present treatment, in that it will make it possible to equivalently describe  any monotone neighbourhood frame only in terms of the minimal neighbourhoods of its states, as detailed in the following paragraph.
\end{remark}

\paragraph*{Finite monotone neighbourhood frames and finite neighbourhood frames.}
\label{para : two sorted frame neighbourhood frame}
Our main focus of interest in the present paper are  finite lattices and their related structures, which are also finite.
For any finite monotone neighbourhood frame
$\mb{F} = (X,\sigma : X \longrightarrow \P\P X)$,
%
%
the collection $\sigma(x)$, which is an upset of $\P X$,
is uniquely identified by the subcollection  of its $\subseteq$-minimal elements.
%
Hence, any $\mb{F}$ as above can be  equivalently represented as the neighbourhood frame $\mb{F}^* := (X,\sigma^*)$ where 
$\sigma^* : X \longrightarrow \P\P X$    maps each state $x$ to   
the $\subseteq$-minimal elements of the collection $\sigma(x)$.
Conversely, any finite neighbourhood frame $\mb{F} =(X,\sigma)$ can be associated with a monotone neighbourhood frame $\mb{F}' := (X, \sigma')$ where $\sigma'(x) = {\uparrow_{\subseteq}} \sigma(x)$ for any $x\in X$,
and moreover, $(\sigma^*)'=\sigma$ for any finite monotone neighbourhood frame.
This correspondence extends to models as follows:
for any finite monotone neighbourhood model $\mb{M} = (\mb{F},v)$, let 
$\mb{M}^* := (\mb{F}^*, v)$ denote its associated finite  neighbourhood model. 
Conversely, for any finite neighbourhood model $\mb{M} = (\mb{F},v)$, let 
$\mb{M}' := (\mb{F}', v)$ denote its associated finite monotone neighbourhood model. 
Thanks to the slightly  non-standard definition of the interpretation of  $\mc{L}_{MML}$-formulas adopted in the present paper (cf.\ Definition \ref{def : semantic MML} and remark \ref{rem : semantic MML}), this 
equivalent representation behaves well with respect to the interpretation of the monotone modal operators. 
Indeed, it is easy to show that for every $\varphi \in \mc{L}_{MML}$, every finite monotone neighbourhood model $\mb{M}$, and every finite neighbourhood model $\mb{N}$,
$$\mb{M},w \Vdash \varphi \quad \tiff \quad \mb{M}^*,w \Vdash \varphi  \quad \quad \quad \quad \tand \quad \quad \quad \quad \mb{N},w \Vdash \varphi \quad \tiff \quad \mb{N}',w \Vdash \varphi  $$
The proof is done by induction on $\varphi$. We do not give it in full, and only report on the case of $\mb{M}$ and the connectives $\ef$ and $\fe$.
        \begin{align*}
        \mb{M},w \Vdash ( \exists \forall ) \varphi
        &\quad\text{ iff }\quad \text{ there exists some }
         C \in \sigma(w) \text{ such that } C \subseteq v( \varphi) \;
         \\
         & \quad \tiff \quad \text{ there exists some }
         C \in min_{\subseteq}\sigma(w) \text{ such that } C \subseteq v( \varphi)
        \\
        & \quad \tiff \quad \mb{M}^{*},w \Vdash ( \exists \forall ) \varphi.
        \\
        ~\\
        \mb{M},w \Vdash ( \forall \exists ) \varphi
        &\quad\text{ iff }\quad
        \text{ for each } C \in \sigma(w)  , \;  C \cap  v( \varphi) \neq \emptyset  
        \\
        &\quad\text{ iff }\quad
        \text{ for each } C \in min_{\subseteq}\sigma(w)  , \;  C \cap  v( \varphi) \neq \emptyset
        \\
        & \quad \tiff \quad \mb{M}^{*},w \Vdash ( \exists \forall ) \varphi  .
        \end{align*}
%
%


\paragraph*{Join-presentations as monotone neighbourhood frames.}
\label{para : join pres and mml frames}

Join-presentations (cf.\ Definition \ref{def : join presentation}) of finite lattices  bear a very close resemblance to neighbourhood frames. This resemblance can be spelled out more precisely, which is what  we are going to do next.

For any finite lattice $L$, let $(J(L),\leq, \M )$ be its join-presentation.
The  monotone neighbourhood frame associated with $L$ is the tuple $\F_L:=(J(L),  \sigma_{\M} : J(L) \longrightarrow \P\P J(L) )$
such that for each $j\in J(L)$,
\begin{equation}
\label{eq : sigma and M}
\sigma_{\M}(j) :=  \{ S \in \P \P J(L) \mid \:  C \subseteq S \text{ for some } C\in \M(j)  \}.
\end{equation}

Clearly, $\sigma_{\M}(j)$ is upward-closed, hence the construction above is well defined. Moreover, since $\M(j)$ is a $\ll$-antichain (see lemma \ref{lem : minimal covers}.1), for all $C$ and $C'$ in $\M(j)$, if $ C \subseteq \:  C'$ then $C=C'$. This immediately implies that $\M(j) $ is the collection $\min_{\subseteq}\sigma_{\M}(j) $  of the $\subseteq$-minimal elements of $\sigma_{\M}(j)$. 

Notice that the construction associating a neighbourhood frame with the join-presentation of a finite lattice $L$, involves a loss of information. Namely, the order $\leq_J$ on the set $J(L)$ of the join-irreducible elements of $L$ cannot be retrieved from the neighbourhood frame $\mb{F}_L$.

For every $L$, we are only interested in valuations on $\mb{F}_L$ which are the dual counterparts of assignments on $L$.
Recall that $L$ is isomorphic to the lattice $\mfLL$
of closed sets of the join-presentation associated with $L$.
Hence, 
we are only interested in valuations
mapping atomic propositions to \textit{closed} subsets, rather than to arbitrary subsets of $\mb{F}_L$.  This motivates the following definition. 
\begin{definition}
\label{def : closed valuation}
For any finite lattice $L$, let a \textit{model} on $\F_L $ be a tuple   $\mb{M}_L =(\F_L , v^*)$ such that $v^*: \AtProp \longrightarrow \mfLL$. We refer to such maps as \textit{closed valuations}.
Then, abusing terminology, the local and global validity of formulas and inequalities on $\F_L$ will be understood relative to closed valuations, that is:
\begin{align*}
\F_L ,j \Vdash \varphi \quad &\tiff \quad (\F_L , v^*),j \Vdash \varphi \text{ for any closed valuation }v^* .\\
\F_L \Vdash \varphi \quad &\tiff \quad \F_L  ,j \Vdash \varphi \text{ for any }j\in J(L).\\
\F_L \Vdash \varphi \leq \psi \quad &\tiff \quad  \text{for any closed valuation } v^* \text{ and any } j\in J(L), \text{ if } (\F_L ,v^*) ,j \Vdash \varphi \text{ then } (\F_L ,v^*) ,j \Vdash \psi.
\end{align*}
\end{definition}

Let us spell out in detail the correspondence between assignments on $L$ and closed valuations on $\F_L$. Clearly, given a set of variables $\AtProp$, closed valuations of $\AtProp$ on $\F_L$ can be identified with assignments of $\AtProp$  on $\mfLL$.
The isomorphism  $\mf{L} : L \longrightarrow \mfLL$ defined by the mapping $a \longmapsto \{ j \in J(L) \mid j\leq a\}$,
with inverse defined by the mapping $S \longmapsto \bigvee_L S$, induce  bijections between assignments on $L$ and assignments on $\mfLL$, defined by post-composition. That is, 
any assignment $v : \AtProp \longrightarrow L$ gives rise to the assignment  $v^* : \AtProp \longrightarrow \mfLL$,
such that  for any $x\in \AtProp$, 
\begin{equation}
\label{eq : def v star}
v^*(x) := \{j\in J(L)\mid j\leq v(x)\}.
\end{equation}  The inverse correspondence maps any assignment/closed valuation $u : \AtProp \longrightarrow \mfLL$ to an assignment $u^{\prime} : \AtProp \longrightarrow L$ such that for any $x\in \AtProp$,  
\begin{equation}
\label{eq : def u prime}
u^{\prime} (x) :=  \bigvee_L u (x).
\end{equation}  
Thus, $v^{* \prime} = v$, and $u^{\prime *} = u$ for any assignment $v$ on $L$ and any assignment $u$ on $\mfLL$.
Hence, for all lattice terms $s$ and $t$ over $\AtProp$,
for any assignment $v$ on $L$ and any assignment $u$ on $\mfLL$,
\begin{align}
    L,v \models s\leq t \quad &\tiff \quad 
    \mfLL, v^* \models s\leq t, \label{equivalence v vstar}\\
    L,u' \models s\leq t \quad &\tiff \quad 
    \mfLL, u \models s\leq t. \label{equivalence u uprime}
\end{align}

\subsection{The standard translation} 
\label{subsec : standard translation}
Thanks to the duality of the previous subsection, and to the correspondence environment introduced above, we are now in a position to define the `standard translation' $ST$ from the language of
lattices to the language of monotone modal logic. The aim of this translation is to have, for any lattice term $t$, any finite lattice $L$, any $j\in J(L)$
and any $v : \AtProp \longrightarrow L$,
\begin{equation}
\label{eq:aim ST}
 L,v\models j\leq t \quad \mbox{ iff }\quad \F_L,v^* ,j\Vdash ST(t),
 \end{equation}
 where $v^* $ is defined as in the discussion after definition \ref{def : closed valuation}.

 The definition of $ST$ pivots on the duality between lattices and join-presentations.
  Namely, any given  interpretation of a lattice term $t$ on a finite lattice $L$ translates to an interpretation of $t$ into the lattice $\mfLL$ of the closed sets of
the join-presentation $(J(L),\leq,\M)$ associated with $L$, via  the fact that $L$
is isomorphic to $\mfLL$. Then, by dually characterizing the interpretation of $t$ in $\mfLL$, we retrieve the interpretation  of  $t$ into the join-presentation $(J(L),\leq,\M)$. In its turn, this interpretation boils down to the satisfaction clause, on $\F_L$, of certain formulas belonging to a \textit{ fragment} of monotone modal logic, which can be recursively defined as follows:
$$ \varphi ::= \bot \mid \top \mid p
\mid \varphi \wedge \varphi
\mid ( \exists \forall ) ( \varphi \vee \varphi ).$$
Let us define $ST$ by the following recursion:
\label{page : def standard translation}
\begin{align*}
ST(p) & =  p\\
ST(\top) & =  \top\\
ST(\bot) & =  \bot\\
ST(t\wedge s) & =  ST(t)\wedge ST(s)\\
ST(t\vee s) & =  \ef (ST(t)\vee ST(s)).\\
\end{align*}
The definition above  recasts \cite[Definition 7.1]{duality} into the language of monotone modal logic.

 In what follows, we will find it useful to expand our propositional language with individual variables of a different sort than propositional variables. These new variables, denoted $\nomj$, $\nomk$, possibly with sub- and superscripts, are to be  interpreted as join-irreducible elements of finite lattices. Let $\Nom$ (for nominals) be the collection of such variables, and   let $\Var := \AtProp \cup \Nom$. Finite lattice assignments from $\Var$ are maps $v : \Var \longrightarrow L$
 such that $v(\nomj) \in J(L)$ for every $\nomj \in \Nom$.
 Each such lattice assignment corresponds to a valuation from $\Var$ to $\F_L$ as described in the discussion at the end of subsection \ref{subsec : MML}. 
\begin{proposition}
\label{prop : local satisfaction}
Let $L$ be a finite lattice which is different from the one-element lattice. Then, for any lattice term $t$ over $\AtProp$, any $j \in J(L)$, and any assignment $v : \Var \longrightarrow L$ with $v(\nomj)=j$,
\begin{equation}
\label{eq:aim ST proposition local}
 L,v\models \nomj\leq t \quad \mbox{ iff }\quad \F_L,v^* ,j\Vdash ST(t),
 \end{equation}
\end{proposition}

\begin{proof}
By induction on $t$. If $t = \top, \bot$, then the statement is clearly true.
If $t=p \in \AtProp$, then $ST(p)=p$.
Then,  
the following chain of logical equivalences holds:
\begin{align*}
L,v\models \nomj\leq p \quad &\tiff \quad v(\nomj) \leq_L v(p) \\
& \tiff \quad \{ k \in J(L) \mid k \leq v(\nomj)  \} \subseteq  \{ k \in J(L) \mid k \leq v(p) \}   \tag{$v(\nomj) \in\{ k \in J(L) \mid k \leq v(\nomj)  \}$ }\\
& \tiff \quad v^*(\nomj) \subseteq v^*(p) \tag{definition of $v^*$} \\
& \tiff \quad v(\nomj) \in v^*(p)  \tag{$v^*(p)$ is a downset}\\
& \tiff \quad \F_L,v^* ,j\Vdash p. \tag{$v(\nomj)=j$}
\end{align*}
The inductive step $t = t_1 \W t_2$
straightforwardly follows from the induction hypothesis.

As for the case $t = t_1 \V t_2$, assume that the equivalence (\ref{eq:aim ST proposition local}) holds for $t_1$ and $t_2$, for every $c \in J(L)$ and for any $v : \Var \longrightarrow L $.
As discussed in the previous subsection (see equation (\ref{equivalence v vstar})), we have
\begin{equation*}
L,v\models \nomj \leq t_1 \V t_2  \quad \tiff \quad \mfLL,v^*\models \nomj \leq t_1 \V t_2. 
\end{equation*}

Let us recall that the meet $\wedge^{*}$ and join $\vee^{*}$ of  $\mfLL$ are  respectively defined as follows:
for all $S, T\in \mfLL$,
\begin{align*}
T \wedge^* S &= T \cap S & &\tand & T \vee^* S & = \overline{T \cup S}
\end{align*}
where $\overline{T \cup S}$ is defined in \eqref{eq : def closure}, that is:
$ \overline{T \cup S} = \{ x\in X \mid \exists C\in \mc{M}(x) \: : \: C \subseteq T \cup S\}. 
$
Hence, the following chain of logical equivalences holds:
\begin{align*}
& \phantom{\tiff \quad}\mfLL,v^*\models \nomj \leq t_1 \V t_2 \\
& \tiff \quad v^*(\nomj) \subseteq v^*(t_1 \V t_2)  \\
& \tiff \quad j = v(\nomj) \in v^*(t_1 \V t_2)  \\
& \tiff \quad j \in\overline{v^*(t_1) \cup v^*(t_2)} \\
& \tiff \quad \text{there exists some } C\in \M(j) \text{ such that } c \in v^*(t_1) \tor c \in v^*(t_2) \text{ for all }c\in C  \\
& \tiff \quad \text{there exists some } C\in \M(j) \text{ such that } {\downarrow}c \subseteq v^*(t_1) \tor {\downarrow}c \subseteq v^*(t_2) \text{ for all }c\in C,
\end{align*}
where ${\downarrow}c := \{ k\in J(L) \mid k\leq c\}$.
For any $c\in J(L)$, let $u_c$ be the $\nomj$-variant of $v^*$ such that $u_c(\nomj) = {\downarrow} c$. Hence, the previous clause can be equivalently rewritten as follows:
\begin{align*}
& \text{there exists some } C\in \M(j) \text{ such that  for all }c\in C,  \quad \mf{L}_L, u_c \models \nomj \leq t_1  \;\tor\; \mf{L}_L, u_c \models \nomj \leq t_2.
\end{align*}
By equation (\ref{equivalence u uprime}), the clause above can equivalently rewritten as follows:
\begin{align*}
& \text{there exists some } C\in \M(j) \text{ such that  for all }c\in C,  \quad L, u_c' \models \nomj \leq t_1  \;\tor\;  L, u_c' \models \nomj \leq t_2. 
\end{align*}
By the induction hypothesis, the clause above is equivalent to
the following one:
\begin{align*}
& \text{there exists some } C\in \M(j) \text{ such that  for all }c\in C,  \quad \F_L, (u_c')^*, u_c'(\nomj) \Vdash  ST(t_1)  \;\tor\;  \F_L, (u_c')^*, u_c'(\nomj) \Vdash  ST(t_2). 
\end{align*}
Moreover, as discussed after definition (\ref{def : closed valuation}), we have that $u_c'(\nomj) = \bigvee_L u_c(\nomj) = \bigvee_L {\downarrow}c = c$, and $(u_c')^*=u_c$. Hence, the clause above can be simplified as follows:
\begin{align*}
& \text{there exists some } C\in \M(j) \text{ such that  for all }c\in C,  \quad \F_L, u_c, c \Vdash  ST(t_1)\;\tor\;  \F_L, u_c, c \Vdash  ST(t_2), 
\end{align*}
 and then as follows:
\begin{align*}
& \text{there exists some } C\in \M(j) \text{ such that  for all }c\in C,  \quad c \in u_c(  ST(t_1))\;\tor\;  c \in u_c(  ST(t_2)). 
\end{align*}
Since $t_1$ and $t_2$ are lattice terms over $\AtProp$, no nominal variable occurs in them, and hence 
$u_c(  ST(t_1)) = v^*(  ST(t_1))$ and $u_c(  ST(t_2)) = v^*(  ST(t_2))$.
Thus, we can equivalently rewrite the clause above as follows:
\begin{align*}
& \text{there exists some } C\in \M(j) \text{ such that  for all }c\in C,  \quad c \in v^*(  ST(t_1))\;\tor\;  c \in v^*(  ST(t_2)).
\end{align*}
By (\ref{eq : sigma and M}), and  since $\M(j) = \min_{\subseteq}\sigma_{\M}(j) $ (see discussion below  (\ref{eq : sigma and M})), the condition above  is equivalent to
\begin{align*}
& \text{there exists some } S\in \sigma_{\M}(j) \text{ such that  for all }c\in S,  \quad c \in v^*(  ST(t_1))\;\tor\;  c \in v^*(  ST(t_2)).
\end{align*}
By definition, this is equivalent to
\begin{align*}
& \F_L,v^* ,j\Vdash \ef(ST(t_1) \V ST(t_2)),
\end{align*}
as required.
\end{proof}

The following corollary gives semantic justification to the standard translation, and provides the mathematical basis for our general approach of
obtaining dual characterization results for  finite lattices 
as instances of correspondence arguments in the language of monotone modal logic.
Recall that, by definition \ref{def : closed valuation}, 
$$\F_L \Vdash \varphi \leq \psi \quad \tiff \quad  \text{for any closed valuation } v^* \text{ and any } j\in J(L), \text{ if } (\F_L ,v^*) ,j \Vdash \varphi \text{ then } (\F_L ,v^*) ,j \Vdash \psi.$$

\begin{corollary}
\label{cor : lattice frame equi}
Let $L$ be a finite lattice. Then, for every lattice term $t$ and $s$,  
\begin{equation*}
\label{eq : ST}L \models t\leq s  \quad
     \tiff \quad \F_L \Vdash ST(t)\leq ST(s).
\end{equation*} 
\end{corollary}

\begin{proof}
Notice that finite lattices are join-generated by their join-irreducible elements.
Hence, the condition $L \models t\leq s$ 
is equivalent to the following:
\begin{equation}
\label{eq : for any valuation stuff 1} 
\text{for any assignment } v: \AtProp \longrightarrow L,\; \and \text{for any } j\in J(L),\;
\tif j \leq v(t) \;\text{ then }j\leq v(s).
\end{equation}
Clause \eqref{eq : for any valuation stuff 1}  is equivalent to the following condition holding 
for any $j\in J(L)$, and 
for any valuation $v: \AtProp \cup \Nom \longrightarrow L$ such that $v(\nomj) = j$:
\begin{equation}
\label{eq : for any valuation stuff 2}
 \text{if } L,v \models \nomj\leq t \text{, then }L,v \models \nomj\leq s.
\end{equation}
By proposition (\ref{prop : local satisfaction}), clause \eqref{eq : for any valuation stuff 2}  is equivalent to: 
\begin{equation}
\label{eq : F v}
\text{if } \F_L, v^*, v(\nomj) \Vdash ST( t) \text{, then } \F_L, v^*, v(\nomj) \Vdash ST( s).
\end{equation}
Next, we claim that clause \eqref{eq : F v} holding 
for any $j\in J(L)$ and 
for any valuation $v: \AtProp \cup \Nom \longrightarrow L$ such that $v(\nomj) = j$ is equivalent to the following: 
\begin{equation*}
\label{eq : F u}
\text{for any }j\in J(L), \text{ and for any closed valuation } u : \AtProp  \longrightarrow \mf{L}_L,\;
\text{if } \F_L, u, j \Vdash ST( t) \text{, then } \F_L, u, j \Vdash ST( s).
\end{equation*}
The latter condition is equivalent to $ \F_L  \Vdash ST(t) \leq ST( s)$, as desired.

To finish the proof, let us prove the claim.
For the direction from top to bottom,  fix a closed valuation $u : \AtProp  \longrightarrow \mf{L}_L$ such that $\F_L, u, j \Vdash ST( t) $ and let $v : \AtProp \cup \Nom \longrightarrow  L$
coincide with $u'$ on $\AtProp$ (cf.\ \eqref{eq : def u prime}) and be such that $v(\nomj)=j$.
By assumption, (\ref{eq : F v}) holds for our choice of $v$. Since $(u')^* = u$, we have that
$v^*$ coincides with $u$ on $\AtProp$, hence
$\F_L, v^*, v(\nomj) \Vdash ST( t) $. Then, by (\ref{eq : F v}), $\F_L, v^*, v(\nomj) \Vdash ST( s) $.
Since $v^*$ coincides with $u$ on $\AtProp$, we have $\F_L, u, j \Vdash ST( s) $ as required.
The direction from bottom to top is proved similarly.
\end{proof}


\subsection{An informal presentation of the algorithm ALBA}
\label{subsec : intro ALBA}


In the present subsection, we illustrate how ALBA works. Our presentation is based on \cite{ALBA, unifiedCorrespondence, CFPS14}. Rather than presenting the algorithm formally, in what follows we will run ALBA on one of the best known examples in correspondence theory, namely $\Diamond \Box p \rightarrow \Box \Diamond p$. It is well known that for every Kripke frame $\mathcal{F} = (W, R)$,
\[
\mathcal{F} \Vdash \Diamond \Box p \rightarrow \Box \Diamond p \,\,\,\mbox{ iff }\,\,\, \mathcal{F} \models \,\forall xyz\,(Rxy \wedge Rxz \rightarrow \exists u(Ryu \wedge Rzu)).
\]
As is discussed at length in \cite{ALBA, unifiedCorrespondence}, every piece of the argument used to prove this correspondence on Kripke frames can be translated by duality  to their complex algebras (cf.\ \cite[Definition 5.21]{Blackburn}), which, as is well known, are complete atomic boolean algebras with operators. We will show how this is done in the case of the example above. First of all, the above validity condition on 
$\mathcal{F}$ translates to its complex algebra $\mathbb{A}$ as $\llb \Diamond\Box p  \rrb \subseteq \llb \Box \Diamond p \rrb $ for every assignment of $p$ into $\mathbb{A}$, so this validity clause can be rephrased as follows:
\begin{equation}\label{Church:Rosser}
\mathbb{A} \models \forall p [\Diamond\Box p \leq \Box\Diamond p].
\end{equation}
Since, in a complete atomic boolean algebra, every element is both the join of the completely join-prime elements (the set of which is denoted $\jty(\mathbb A)$) below it and the meet of the completely meet-prime elements (the set of which is denoted $\mty(\mathbb A)$) above it, the condition above can be equivalently rewritten as follows:
%
\begin{equation*}
\mathbb{A} \models \forall p [\bigvee\{i\in \jty(\mathbb{A})\mid i \leq\Box\Diamond p\} \leq \bigwedge\{m\in \mty(\mathbb{A})\mid\Box\Diamond p\leq m\}].
\end{equation*}
By elementary properties of least upper bounds and greatest lower bounds in posets (cf.\ \cite{DaveyPriestley2002}), this condition is true if and only if every element in the join is less than or equal to every element in the meet. Thus, the condition above  can be equivalently rewritten as:
\[
\mathbb{A} \models \forall p \forall \nomi \forall \cnomm [(\nomi \leq\Diamond \Box p \,\,\,\& \,\,\,\Box \Diamond p\leq \cnomm) \Rightarrow \nomi \leq \cnomm ],
\]
where the variables $\nomi$ and $\cnomm$ range over $\jty(\mathbb{A})$ and $\mty(\mathbb{A})$ respectively. 
Since this presentation is geared towards the treatment in section \ref{sec : Dchain}, we find it useful to slightly depart from the standard treatment in \cite{ALBA} and eliminate the conominal $\cnomm$ as follows. First, notice that the clause above is clearly equivalent to the following clause:
\[
\mathbb{A} \models  [\exists p \exists \nomi \exists \cnomm (\nomi \leq\Diamond \Box p \,\,\,\& \,\,\,\Box \Diamond p\leq \cnomm
\,\,\,\& \,\,\, \nomi \nleq \cnomm ) ]\Rightarrow \false .
\]
Second, notice that, in any complete atomic boolean algebra $\mathbb{A}$, for each $i\in \jty(\mathbb{A}) $ and each $m \in \mty (\mathbb{A})$, one has $i\nleq m$ iff $m = \kappa(i)$, where $\kappa(i) = \bigvee \{ j \in \jty(\mathbb{A}) \mid j \neq i   \} \in \mty(\mathbb{A})$.  
Hence, the clause above is equivalent to the following clause:
\begin{equation}
\label{eq : clause kappa(i)}
\mathbb{A} \models  [\exists p \exists \nomi  (\nomi \leq\Diamond \Box p \,\,\,\& \,\,\,\Box \Diamond p\leq \kappa(\nomi) ) ]\Rightarrow \false .
\end{equation}
Since $\mathbb{A}$ is in particular atomistic, the element of $\mathbb{A}$ interpreting $\Box p$ is the join of the completely join-prime elements below it. Hence, if $i \in \jty(\mathbb A)$ and $i \leq \Diamond \Box p$, because $\Diamond$ is completely join-preserving on $\mathbb{A}$, we have that
\[
i \leq \Diamond(\bigvee \{j \in \jty(\mathbb{A})\mid j\leq  \Box p \}) = \bigvee\{\Diamond j \mid j\in \jty(\mathbb{A})\mbox{ and } j\leq \Box p \},
\]
which implies that $i \leq \Diamond j_0$ for some $j_0 \in \jty(\mathbb A)$ such that $j_0 \leq \Box p$. Hence, we can equivalently rewrite the validity clause (\ref{eq : clause kappa(i)}) as follows:
\[
\mathbb{A} \models   [\exists p \exists \nomi (\exists \nomj(\nomi\leq\Diamond \nomj \,\,\,\& \,\,\, \nomj\leq \Box p) \,\,\,\& \,\,\, \Box\Diamond p \leq \kappa(\nomi)) ]\Rightarrow \false ,
\]
and then as follows:
\[
\mathbb{A} \models \forall p  \forall \nomi \forall \nomj [ (\nomi\leq\Diamond \nomj \,\,\,\& \,\,\, \nomj\leq \Box p \,\,\,\& \,\,\,\Box\Diamond p \leq \kappa(\nomi) )\Rightarrow \false ].
\]
Now we observe that the operation  $\Box$ preserves arbitrary meets in $\mathbb{A}$, which is in particular a complete lattice. By the general theory of adjunction in complete lattices, this is equivalent to $\Box$ being a right adjoint (cf.\  \cite[proposition 7.34]{DaveyPriestley2002}). It is also well known that the left adjoint of $\Box$ is the operation $\Diamondblack$, which can be thought of as the backward looking diamond of tense logic. Hence the condition above can be equivalently rewritten as:
\[
\mathbb{A} \models \forall p \forall \nomi  \forall \nomj[(\nomi\leq\Diamond \nomj \,\,\,\& \,\,\, \Diamondblack \nomj \leq p  \,\,\,\& \,\,\, \Box\Diamond p \leq \kappa(\nomi)) \Rightarrow \false ],
\]
and then as follows:
\begin{equation}
\label{Before:Ack:Eq}
\mathbb{A} \models \forall \nomi  \forall \nomj[(\nomi\leq\Diamond \nomj \,\,\,\& \,\,\, \exists p (\Diamondblack \nomj \leq p  \,\,\,\& \,\,\, \Box\Diamond p \leq \kappa(\nomi)) )\Rightarrow \false ].
\end{equation}
At this point we are in a position to eliminate the variable $p$ and equivalently rewrite the previous condition as follows:
\begin{equation}
\label{After:Ack:Eq}
\mathbb{A} \models \forall \nomi  \forall \nomj[(\nomi\leq\Diamond \nomj \,\,\,\& \,\,\, \Box\Diamond \Diamondblack \nomj \leq \kappa(\nomi)) \Rightarrow \false ].
\end{equation}

Let us justify this equivalence: for the direction from top to bottom, fix an interpretation $v$, and assume that $\mathbb{A},v \models \nomi\leq\Diamond \nomj$ and $\mathbb{A},v \models  \Box \Diamond \Diamondblack \nomj \leq \kappa(\nomi)$. Consider the $p$-variant $v^\ast$ of $v$ such that $v^\ast(p) = \Diamondblack \nomj$. Then it can be easily verified that $\mathbb{A}, v^{\ast} \models \nomi\leq\Diamond\ \nomj $ and $\mathbb{A}, v^{\ast} \models  \Diamondblack \nomj \leq p $ and $\mathbb{A}, v^{\ast} \models   \Box \Diamond p \leq \kappa(\nomi))$,
which by assumption leads to an inconsistency.

Conversely,  fix an interpretation $v$ such that
$\mathbb{A},v \models \nomi\leq\Diamond \nomj $ and $\mathbb{A},v \models\exists p (\Diamondblack \nomj \leq p  \,\,\,\& \,\,\, \Box \Diamond p \leq \kappa(\nomi)) $. Then, by monotonicity,  the antecedent of (\ref{After:Ack:Eq}) holds under $v$, which leads again to an inconsistency. This is an instance of the following result, known as \textit{Ackermann's lemma} (\cite{Ackermann:Untersuchung}, see also \cite{Conradie:et:al:SQEMAI}):
\begin{lemma}
    \label{Right:Ackermann}
    Let $\alpha, \beta(p), \gamma(p)$ be $L$-formulas, such that $\alpha$ is $p$-free, $\beta$ is positive and $\gamma$ is negative in $p$. For any assignment $v$ on an $L$-algebra $\mathbb A$, the following are equivalent:
    \begin{enumerate}
    \item  $\mathbb{A}, v \models \beta(\alpha/p) \leq \gamma(\alpha/p)$;
    \item  there exists a $p$-variant $v^\ast$ of $v$ such that $\mathbb{A}, v^\ast \models \alpha \leq p$ and $\ \mathbb{A}, v^\ast \models \beta(p)\leq\gamma(p)$.
    \end{enumerate}
\end{lemma}
The proof is  similar to that of \cite[lemma 4.2]{ALBA}. Whenever, in a reduction, we reach a shape in which the lemma above (or its order-dual) can be applied, we say that the condition is in \textit{Ackermann shape}\label{Ackermann:Shape}. 

By the definition of $\kappa(\nomi)$, the inequality $\Diamond \Box \Diamondblack \nomj \leq \kappa(\nomi))$ is equivalent to 
$\nomi \nleq \Diamond \Box \Diamondblack \nomj$.
 Hence,  clause (\ref{After:Ack:Eq}) can be equivalently rewritten as follows:
\begin{equation}
\mathbb{A} \models \forall \nomi  \forall \nomj[(\nomi\leq\Diamond \nomj \,\,\,\& \,\,\, \nomi \nleq  \Box \Diamond \Diamondblack \nomj ) \Rightarrow \false ],
\end{equation}
and then as follows:
\begin{equation}
\mathbb{A} \models \forall \nomi  \forall \nomj[\nomi\leq\Diamond \nomj   \Rightarrow \nomi \leq  \Box \Diamond \Diamondblack \nomj  ].
\end{equation}
By the atomicity of $\mathbb{A}$, the clause above is equivalent to:
\begin{equation}
\mathbb{A} \models  \forall \nomj[\Diamond \nomj \leq   \Box \Diamond \Diamondblack \nomj  ].
\end{equation}
By again applying the fact that $\Box$ is a right adjoint we obtain

\begin{equation}
\label{Eq : reduced}
\mathbb{A} \models \forall \nomj[\Diamondblack \Diamond \nomj \leq \Diamond \Diamondblack \nomj].
\end{equation}
%
%
Recalling that $\mathbb{A}$ is the complex algebra of $\mathcal{F} = (W,R)$, this gives $\forall w( R[R^{-1}[ w ]] \subseteq  R^{-1}[R[ w ]]$. Notice that  $R[R^{-1}[ w ]]$ is the set of all states $x \in W$ which have a predecessor $z$ in common with $w$, while  $R^{-1}[R[ w ]]$ is the set of all states $x \in W$ which have a successor in common with $w$. This can be spelled out as
\[
\forall x \forall w ( \exists z (Rzx \wedge Rzw) \rightarrow \exists y (Rxy \wedge Rwy))
\]
or, equivalently,
\[
\forall z \forall x \forall w ( (Rzx \wedge Rzw) \rightarrow \exists y (Rxy \wedge Rwy))
\]
which is the familiar Church-Rosser condition.



\section{Algorithmic correspondence for monotone modal logic}
\label{sec : 2sorted}
    
A key intermediate step of the present paper is to adapt the algorithm or calculus for correspondence ALBA to monotone modal logic. The interest of this adaptation is independent from the applications to the theory of finite lattices. So, for the sake of modularity and generality,
we work in a more abstract setting than the one associated with finite lattices, to which this adaptation will be applied. 
The general strategy underlying this adaptation is to exploit the well known fact  that  the `exists/for all' and `for all/exists' quantification patterns in the standard interpretation of the monotone modal operators make it possible to regard monotone modal operators as suitable concatenations of {\it normal} modalities. 
This same observation inspired Helle Hansen's syntactic translation \cite[Definition 5.7]{HelleThesis} on which her Sahlqvist correspondence theorem for monotone modal logic is based.
 The present section  is aimed at making all this precise.
In the next subsection, we introduce two-sorted frames, their associated normal modal language, and first order correspondence language. We also spell out the relationship between two-sorted frames and monotone neighbourhood frames, which allows to interpret monotone modal logic on two-sorted frames.
In subsection \ref{subsec : ALBA on two-sorted frames}, we introduce the basic adaptation of ALBA to the normal modal language of two-sorted frames.
 
\subsection{Two-sorted frames}
\label{subsec : two-sorted frames}
 
 \begin{definition}
 \label{def : two sorted frame}
A  \textit{two-sorted frame}  is a structure $\mb{X} = \langle X, Y, \rxy, \ryx\rangle$ such that $X$ and $Y$ are sets,  $\rxy \subseteq X\times Y$
    and $\ryx \subseteq Y\times X$.
    \end{definition}

The existence of the equivalent representation of any finite monotone neighbourhood frame $\mb{F}$ in terms of the finite neighbourhood frame $\mb{F}^*$ (cf.\ paragraph page \pageref{para : two sorted frame neighbourhood frame})
implies that we can equivalently encode any monotone neighbourhood frame $\mb{F}$ as the following two-sorted structure
$( X , Y ,  \rxy , \ryx) $, such that  $Y = \P X$, and  for every $x\in X$ and $y\in Y$,
$$x\rxy y \; \tiff \; y\in  \min_{\subseteq} \sigma(x)
\quad\quad \tand \quad \quad
y\ryx x \; \tiff \; x\in y.$$
The definitions above imply  that $\rxy [x] = \min_{\subseteq} \sigma(x)$  for any $x\in X$, and $\ryx [y] = y$ for any $y\in Y$. 
In the remainder of the paper, for any relation $S \subseteq X\times Y$,
we sometimes use the symbols $x S$ and $S \hspace*{-0.5mm} y$ to denote the sets 
$S[x]$ and $S^{-1}[y]$
respectively.


\label{page : relations and their modalities}
As is well known, each of the two relations $\rxy$ and $\ryx$ gives rise to a pair of semantic normal modal operators:

\begin{center}
\begin{tabular}{r c l r c l}
$\langle \rxy \rangle :  \P Y$ & $\longrightarrow$ & $\P X$ & $[ \rxy ] :  \P Y $ & $\longrightarrow $ &$ \P X$\\
$T $ & $\longmapsto $ & $\rxy^{-1}[T]$ & $T $ & $\longmapsto $ & $(\rxy^{-1}[T^c])^c$\\
 &&&&&\\
$\langle\ryx \rangle  :   \P X $ & $\longrightarrow $&$ \P Y$ & $[ \ryx ] :  \P X $ &$\longrightarrow $ &$\P Y$\\
$S $ &$\longmapsto $ &$\ryx^{-1}[S]$ & $S$ &$ \longmapsto $ &$ (\ryx^{-1}[S^c])^c$ \\
\end{tabular}
\end{center}
where
\begin{center}
\begin{tabular}{c c }
$\rxy^{-1}[T]:=\{  x \in X \mid x\rxy \cap T\neq\varnothing \}$ & $(\rxy^{-1}[T^c])^c:=\{  x \in X \mid x\rxy \subseteq T \}$ \\
$\ryx^{-1}[S]:=\{  y \in y \mid y\ryx \cap S\neq\varnothing \}$ & $(\ryx^{-1}[S^c])^c:=\{  x \in X \mid x\ryx \subseteq S \}$.
\end{tabular}
\end{center}

\begin{definition}
The {\it complex algebra} of the two-sorted frame $\mb{X}$ as above is the tuple
\[(\P X, \P Y, \langle \rxy \rangle, [\rxy], \langle \ryx \rangle, [\ryx]).\]
To make definitions and calculations more readable
we introduce the following convention:
we note  $\leqx$ the order on $\P X$
and $\leqy$ the order on $\P Y$.
\end{definition}

Two-sorted frames and their complex algebras will be used as (nonstandard) models for the modal language $\mc{L}_{MML}$ over $\AtProp$ (cf.\ Definition \ref{def : language MML}), the definition of which we report here for the reader's convenience:
\[
\phi \, ::= \, \bot \mid \top \mid p \mid \neg \phi
\mid \phi \vee \phi \mid \phi \wedge \phi
\mid (\exists \forall) \phi \mid (\forall \exists) \phi.
\]
\begin{definition}
\label{def : two sorted model}
A \textit{two-sorted model} is a tuple $M = (\mb{X}, v)$ such that $\mb{X}$ is  a two-sorted frame, and  $v$ is a map $\AtProp \longrightarrow \P X$.
\end{definition}
Given a valuation $v$, its associated extension function is defined by induction as follows:
\begin{center}
\begin{tabular}{r c l c r c l}
    $\llb \bot \rrb_{v,X}$ & $=$ & $ \emptyset$ &$\quad$& && \\

    $\llb \top \rrb_{v,X}$  & $=$ & $X$ &$\quad$& &&
    \\
    $\llb p \rrb_{v,X}$ & $=$ & $ v(p)$&$\quad$& &&
    \\
    $\llb \neg \phi \rrb_{v,X}$ & $=$ & $\llb \phi \rrb ^{c}$&$\quad$& &&
    \\
    $\llb \phi \vee \psi \rrb_{v,X}$ & $=$ & $
    \llb \phi \rrb_{v,X} \cup \llb \psi \rrb_{v,X}$&$\quad$& &&
    \\
    $\llb \phi \wedge \psi \rrb_{v,X}$ & $=$ & $
    \llb \phi \rrb_{v,X} \cap \llb \psi \rrb_{v,X}$&$\quad$
    \\
    $\llb (\exists\forall) \phi \rrb_{v,X}$ & $=$ & $\langle \rxy \rangle \lb \ryx \rb \llb \phi \rrb_{v,X}$& $\quad$ & $(\ast)$
    \\
    $\llb (\forall\exists) \phi \rrb_{v,X}$ & $=$ & $ [\rxy]\langle \ryx \rangle  \llb \phi \rrb_{v,X}$ & $\quad$ & $(\ast\ast)$
    \\
    \end{tabular}
    \end{center}


\subsection{Basic ALBA on two-sorted frames}
\label{subsec : ALBA on two-sorted frames}

In order to adapt ALBA to the setting of two-sorted frames, we need to define the symbolic language
which ALBA will manipulate. Analogously to what has been done in \cite{ALBA}, let us introduce the language $\mc{L}^+$ as follows:
\begin{align*}
\varphi ::=& \; \bot \mid \top \mid p \mid \nomj \mid \cnomm \mid
\unomj \mid \ucnomm \mid \neg \varphi \mid 
\varphi \V \varphi \mid \varphi \W \varphi \mid
\varphi \smallsetminus \varphi \mid 
\varphi \rightarrow \varphi \mid \\
& \; 
\langle \rxy\rangle \varphi \mid [\rxy]
\varphi \mid \langle \ryx\rangle
\varphi \mid [\ryx] \varphi \mid
[\rxy^{-1}] \varphi \mid
\langle \rxy^{-1}\rangle \varphi \mid
 [\ryx^{-1}] \varphi \mid \langle \ryx^{-1}\rangle \varphi,
\end{align*}
where $p\in \AtProp$, $\nomj \in \Nom_X$, $\unomj \in \Nom_Y$, $\cnomm \in \CNom_X$, $\ucnomm \in \CNom_Y$. 
The language above is shaped on the complex algebra of two-sorted frames. In particular, the variables in $\Nom_X$ and $\Nom_Y$ are to be interpreted as atoms of $\P X$ and $\P Y$ respectively, and the variables in $\CNom_X$ and $\CNom_Y$ are to be interpreted as coatoms of $\P X$ and $\P Y$.
Moreover, the interpretation of the modal operators is the natural one suggested by the notation and indeed we are using the same symbols to denote both the operators and their interpretations.
Finally,  clauses $(\ast)$ and $(\ast\ast)$ in Definition \ref{def : two sorted model} justifies the definition
of the obvious translation from formulas of $\mc{L}_{MML}$
to formulas in $\mc{L}^+$.
In what follows, we introduce the ALBA rules which are sound on general two-sorted structures.


\paragraph*{Adjunction and residuation rules.}
\label{subsec : adjunction}

It is well known that, in the setting of boolean algebras,
the interpretation of the conjunction $\wedge$ has a right residual, which is the interpretation of the implication, $\rightarrow$, and
the interpretation of the disjunction $\vee$ has a left residual, which is the interpretation of the subtraction $\smallsetminus$. Thus, the following rules are sound and invertible in the two boolean algebras associated with any two-sorted structure:
$$\dfrac{\alpha \wedge \beta \leq \gamma}{\alpha \leq \beta \rightarrow \gamma} \ajW
\quad\quad\quad\quad
\dfrac{\alpha \leq \beta \vee \gamma}{\alpha \smallsetminus \beta \leq \gamma} \ajV$$
Moreover, it follows  from very  well known  facts in modal logic that, 
for any two-sorted structure,   $\Drxy$ 
(resp.\ $\Brxy$) has a right (resp.\ left) adjoint, which is $\Bryx$ (resp.\ $\Dryx$). Thus, the following rules are sound and invertible on any two-sorted structure:
$$\dfrac{\Drxy \alpha \leqx \beta}{\alpha \leqy \Brxym \beta} \tiny{\aj \Drxy}
\quad\quad\quad\quad
\dfrac{\alpha \leqx \Brxy \beta}{\Drxym \alpha \leqy \beta} \tiny{\aj \Brxy}
\quad\quad\quad\quad
\dfrac{\Dryx \alpha \leqy \beta}{\alpha \leqx \Bryxm \beta} \tiny{\aj \Dryx}
\quad\quad\quad\quad
\dfrac{\alpha \leqy \Bryx \beta}{\Dryxm \alpha \leqx \beta} \tiny{\aj} \Bryx$$

\paragraph*{Approximation rules.}
\label{subsec : Approximation}
The soundness and invertibility of the  rules below straightforwardly follows from the complete join- (resp.\ meet-)preservation properties of the modalities $\Drxy, \Brxy, \Dryx$ and $\Bryx$, and also from the fact that the boolean algebras $\P X$ and $\P Y$ are both completely join-generated by their completely join-irreducible elements and  completely meet-generated by their completely meet-irreducible elements.
For more details on this the reader is referred to \cite{ALBA}.
\begin{align*}
&\dfrac{\nomi \leqx \Drxy \alpha}{\exists \ud{\nomj} \; (\nomi \leqx \Drxy \ud{\nomj}  \;\;\&\;\;  \ud{\nomj} \leqy \alpha)} \ap\Drxy
&
&\dfrac{ \Brxy \alpha \leqx \cnomm}{\exists \ud{\cnomn} \; (\alpha \leqy \ud{\cnomn} \;\;\&\;\; \Brxy \ud{\cnomn} \leqx \cnomm)} \ap\Brxy
\\
&\dfrac{\ud{\nomi} \leqy \Dryx \alpha}{\exists \nomj \; (\ud{\nomi} \leqy \Dryx \nomj \;\;\&\;\;  \nomj \leqx \alpha)} \ap\Dryx
&
&\dfrac{ \Bryx \alpha \leqy \ud{\cnomm}}{\exists \cnomn \; (\alpha \leqx \cnomn \;\;\&\;\; \Bryx \cnomn \leqy \ud{\cnomm})} \ap\Bryx
\end{align*}

\paragraph*{Splitting rules.}

The following rules reflect the fact that the logical conjuction and disjunction are respectively interpreted with the greatest lower bound and least upper bound lattice operations, and hence are sound and invertible.

$$\frac{\varphi \leq \psi_1 \W \psi_2}{\varphi \leq \psi_1 
\;\; \& \;\; \varphi \leq \psi_2}  \spW
\quad\quad\quad\quad
\frac{ \psi_1 \V \psi_2 \leq \varphi}{\psi_1 \leq \varphi
\;\; \& \;\;  \psi_2 \leq \varphi}  \spV
$$

\paragraph*{Ackermann rules.}The  soundness and invertibility of the following rules (here below is the right-Ackermann rule)
 has been discussed in \cite[Lemmas 4.2 and 4.3]{ALBA}.

$$\frac{ \exists p \; \Big [ 
\mbox{\Large{\&}}_{i=1}^{n} \{ \alpha_i \leq p   \} \;\; \& \;\;
 \mbox{\Large{\&}}_{j=1}^m \{ \beta_j(p)\leq \gamma_j(p)\} \Big ]
}{ 
\mbox{\Large{\&}}_{j=1}^m 
\{\beta_j(\bigvee_{i=1}^n \alpha_i)\leq \gamma_j(\bigvee_{i=1}^n \alpha_i) \} } (RAR)$$
where
$p$ does not occur in $\alpha_1, \ldots, \alpha_n$, the formulas
$\beta_{1}(p), \ldots, \beta_{m}(p)$ are positive in $p$, and
$\gamma_{1}(p), \ldots, \gamma_{m}(p)$ are negative in $p$.
Here below is the left-Ackermann rule:

$$\frac{\exists p \; \Big [ 
\mbox{\Large{\&}}_{i=1}^{n} \{  p \leq \alpha_i   \} \;\; \& \;\;
 \mbox{\Large{\&}}_{j=1}^m \{ \beta_j(p)\leq \gamma_j(p)\} \Big ]
 }{
 \mbox{\Large{\&}}_{j=1}^m \{ \beta_j(\bigwedge_{i=1}^n \alpha_i)\leq \gamma_j(\bigwedge_{i=1}^n \alpha_i)  \}}
(LAR)$$
where
$p$ does not occur in $\alpha_1, \ldots, \alpha_n$, the formulas
$\beta_{1}(p), \ldots, \beta_{m}(p)$ are negative in $p$, and
$\gamma_{1}(p), \ldots, \gamma_{m}(p)$ are positive in $p$.

\paragraph*{Boolean tautologies.}Clearly, the appropriate boolean and lattice tautologies justify the soundness and invertibility of the following rules.
For the sake of conciseness, some of these rules will be given as formula-rewriting rules rather than as equivalences between inequalities.

\begin{align*}
\frac{\varphi\V \bot}{\varphi} \Vbot
&&
\frac{\varphi \V (\psi_1 \W \psi_2) }{(\varphi \V \psi_1) \W
(\varphi \V \psi_2)}
\distVW
&&
\frac{\neg\neg \phi}{\phi}\tNN
&&
\frac{A \leq B }{(A \W B)  =  A}
\baW
\\
~\\
\frac{\varphi\W \top}{\varphi} \Wtop
&&
\frac{\varphi \W (\psi_1 \V \psi_2) }{(\varphi \W \psi_1) \V
(\varphi \W \psi_2)}
\distWV
&&
\frac{x \V (y \smallsetminus x)}{x \V y} \tiny{\tV} 
&&
\frac{B \leq A}{(A \V B)  =  A}
\baV
\\
~\\
\frac{\varphi \V \psi}{ \psi \V \varphi}
\cmV
&&
\frac{(\varphi \W \psi)\W\chi}{ \psi \W (\varphi \W \chi)}
\asW
&&
\frac{x \W (x \rightarrow y)}{x \W y} \tW
&&
\frac{\varphi \W \psi \leq \bot}{\varphi \leq \neg \psi}\tWB 
\\
~\\
\frac{\varphi \W \psi}{ \psi \W \varphi}
\cmW
&&
\frac{(\varphi \V \psi)\V \chi }{ \psi \V (\varphi\V \chi)}
\asV
&&
\frac{\xi \W (\varphi \smallsetminus \psi ) \leq \chi}{\xi \W
\varphi \leq \psi \V \chi } \tminus
&& 
\frac{\neg (\varphi \V \psi)}{\neg \varphi \W \neg \psi}\DeMorgan
\end{align*}

\paragraph*{Behaviour of atoms.}In any complete atomic boolean algebra,
$\kappa(j)$ coincides with $\neg j$ for each completely join-irreducible element $j$.
Thus, the following rules are sound and invertible in the two boolean algebras associated with any two-sorted structure:
\begin{align*}
\frac{\nomj\W s \leqx\bot}{s\leqx \kappa(\nomj)} \AtCoat
&& \frac{\nomj\W s \leqx \kappa(\nomj)}{s\leqx \kappa(\nomj)} \AtCoattwo
&& \frac{\ud{\nomj}\W s \leqy\bot}{s\leqy \kappa(\ud{\nomj})} \AtCoat
&& \frac{\ud{\nomj} \W s \leqy \kappa(\ud{\nomj})}{s\leqy \kappa(\ud{\nomj})} \AtCoattwo
\\
\frac{\nomj \leq s \V t \quad s \leq \kappa(\nomj)}{\nomj\leq t \quad s \leq \kappa(\nomj)} \MT &&
\frac{\ud{\nomj} \leqy s \V t \quad s \leqy \kappa(\ud{\nomj})}{\ud{\nomj} \leqy t \quad s \leqy \kappa(\ud{\nomj})} \MT 
&&&&
\end{align*}

\paragraph*{Logical rules.}Finally, we find it useful to stress that ALBA is able to perform elementary equivalent simplifications such as those represented in the rules below:

\begin{align*}
\frac{\phi \leq \psi \quad\quad \phi \leq \psi}{\phi \leq \psi} \bis
&&\frac{A=B \quad\quad t(A)\leq s(A)}{A=B \quad\quad t(B)\leq s(B)}\Substitution
&&
\frac{\varphi \leq \psi \quad \psi \leq \chi}{
\varphi \leq \psi \quad \psi \leq \chi \quad \varphi \leq \chi}\tra
\end{align*}
where $t(B)$ and $s(B)$ are obtained by replacing occurrences of $A$ with $B$ in $t$ and $s$ respectively.

\paragraph*{Rules for normal modalities.}
\label{rules for normal modalities} The soundness and invertibility of the following rules $\tBD$, $\tDB$ and $\tnm$ straightforwardly follows from well known validities for classical normal modal logic. The soundness and invertibility of $\tRRminus$ immediately follows from the definition of the semantics of $\langle R \rangle $ and $\langle R^{-1} \rangle $.
$$\frac{\lb R \rb X }{ \neg \langle R \rangle \neg X} \tBD
\quad \quad
\frac{\neg \langle R \rangle \neg X}{ \lb R \rb X } \tDB
\quad\quad 
\frac{\nomj \leq  \langle R \rangle \unomj}{\unomj \leq \langle R^{-1} \rangle \nomj} \tRRminus
\quad\quad 
\frac{X \leq \langle R \rangle A  \quad X \leq \lb R \rb B}{
X \leq \langle R \rangle (A\W B)  \quad X \leq \lb R \rb B }\tnm .
$$

\section{Enhancing the algorithm for correspondence}
\label{sec : Enhancing ALBA}

We are working towards being able to account for Nation's characterisation in \cite{NationApproach} as an instance of algorithmic correspondence for the monotone modal logic language defined in Definition \ref{def : language MML}.
As we saw in subsection \ref{subsec : standard translation},  the
validity of a lattice inequality on any finite lattice $L$ corresponds to 
the validity of the standard translation (cf.\ page \ref{page : def standard translation}) of the given inequality on the join-presentation $\F_L$  associated with $L$ \textit{restricted}  to closed valuations (cf.\ Definition \ref{def : closed valuation}).
However, the  version of ALBA for monotone neighbourhood frames 
defined in the previous section is not equipped to recognize closed valuations and properly treat them.
Therefore, in the present section, we enhance the environment of two-sorted frames with an extra relation which encodes the order on the join-presentation $\F_L$.   
On this environment, additional ALBA rules can be shown to be sound, thanks to which closed valuations can be accounted for.

\subsection{Enriched two-sorted frames}
\label{subsec : Enriched two sorted frame}
In the present subsection, we introduce the enriched two-sorted frames, and we show that the join-presentation of any finite lattice can be equivalently represented as an enriched two-sorted frame.
 \begin{definition}
An  \textit{enriched two-sorted frame}  is a structure $\mb{E} = \langle X, Y, \rxy, \ryx, \rxx \rangle$ such that $\langle X, Y, \rxy, \ryx \rangle$ is a two-sorted frame (cf.\ Definition \ref{def : two sorted frame}),  and $\rxx \subseteq X\times X$. An enriched two-sorted frame is
\begin{itemize}
\item \textit{ordered} if $\rxx$ is a partial order;
\item \textit{minimal} if 
\begin{itemize}
\item it is ordered,
\item $x\rxy y $ implies that the set $y\ryx = \{ x'\in X \mid y\ryx x'\}$ is a $\rxx$-antichain for every $x \in X$ and $y\in Y$,
\item the collection $\{ y\ryx \mid y \in x\rxy \} $ is a $\ll$-antichain for any $x\in X$,
\end{itemize}
where $\ll$ is the refinement relation associated with the partial order $(X,\rxx)$ (cf.\ \eqref{eq : def refinement}).
\item \textit{monotone} if for all $x,x' \in X $, and for each $y\in Y$,
if $x'\rxx x$ and $x \rxy y$, then $y' \ryx  \ll y \ryx $ for  some $y' \in x'\rxy $; 
\item \textit{reflexive} if for every $x\in X$ there exists  some $y\in Y$ such that $x \rxy y $ and $y \ryx \ll \{ x\}$;
\item \textit{transitive} if for every $x\in X$ and $y\in Y$, if $y \in x\rxy $ and   $y_{x'} \in  x' \rxy   $ for some $ x' \in y \ryx $ then there exists some $y'\in Y$ such that $x\rxy y'$ and  $y' \ryx \ll \bigcup \{ y_{x'}\ryx \mid x' \in y\ryx\}$;
\item \textit{direct} if it is ordered, minimal, monotone, reflexive and transitive;
\end{itemize}
    \end{definition}

\begin{definition}
\label{def : two sorted frame of a lattice}
Any join-presentation $\A_L :=(J(L),\leq_J, \M : J(L) \longrightarrow \P\P J(L))$ can be equivalently represented as an enriched two-sorted frame $\mb{E}_L := (X,Y,\rxy,\ryx,\rxx)$
by setting 
\begin{align*}
& X := J(L) ,
\quad 
Y:= \{ S \in \P J(L) \mid S \text{ is a } \leq_J \text{-antichain}  \} ,
\\
& \rxy := \{ (x,y)\in X\times Y \mid y \in \M(x) \},
\quad 
\ryx :=\; \ni,
\quad \tand \quad 
\rxx = \;\leq_J .
\end{align*}

It can be easily verified that for every finite lattice $L$, the enriched two-sorted frame $\mb{E}_L$  is direct.
\end{definition}



Similarly to what has been discussed at the beginning of subsection \ref{subsec : two-sorted frames} (cf.\ page \pageref{page : relations and their modalities}), 
each of the three relations $\rxy$,  $\ryx$, and $\rxx$ gives rise to a pair of semantic normal modal operators:

\begin{center}
\begin{tabular}{r c l r c l}
$\langle \rxy \rangle :  \P Y$ & $\longrightarrow$ & $\P X$ & $[ \rxy ] :  \P Y $ & $\longrightarrow $ &$ \P X$\\
$T $ & $\longmapsto $ & $\rxy^{-1}[T]$ & $T $ & $\longmapsto $ & $(\rxy^{-1}[T^c])^c$\\
 &&&&&\\
$\langle\ryx \rangle  :   \P X $ & $\longrightarrow $&$ \P Y$ & $[ \ryx ] :  \P X $ &$\longrightarrow $ &$\P Y$\\
$S $ &$\longmapsto $ &$\ryx^{-1}[S]$ & $S$ &$ \longmapsto $ &$ (\ryx^{-1}[S^c])^c$ \\
 &&&&&\\
$\langle \rxx \rangle  :   \P X $ & $\longrightarrow $&$ \P X$ & $[ \rxx ] :  \P X $ &$\longrightarrow $ &$\P X$\\
$S $ &$\longmapsto $ &$\rxx^{-1}[S] $ & $S$ &$ \longmapsto $ &$ (\rxx^{-1}[S^c])^c$ \\
\end{tabular}
\end{center}
where
\begin{center}
\begin{tabular}{c c }
$\rxy^{-1}[T]:=\{  x \in X \mid x\rxy \cap T\neq\varnothing \}$ & $(\rxy^{-1}[T^c])^c:=\{  x \in X \mid x\rxy \subseteq T \}$ \\
$\ryx^{-1}[S]:=\{  y \in Y \mid y\ryx \cap S\neq\varnothing \}$ & $(\ryx^{-1}[S^c])^c:=\{  y \in Y \mid y\ryx \subseteq S \}$\\
$\rxx^{-1}[S]  := \{ x\in X \mid  x \rxx \cap S \neq \varnothing  \} $   &
$(\rxx^{-1}[S^c])^c := \{  x \in X \mid x \rxx  \subseteq T \}.
$
\end{tabular}
\end{center}

\begin{definition}
The {\it complex algebra} of the enriched two-sorted frame $\mb{E}$ as above is the tuple
\[(\P X, \P Y, \langle \rxy \rangle, [\rxy], \langle \ryx \rangle, [\ryx], \langle \rxx \rangle, [\rxx]).\]

\end{definition}

\begin{definition}
An \textit{enriched two-sorted model} is a tuple $\mb{M} = (\mb{E}, v)$ such that $\mb{X}$ is  an enriched two-sorted frame, and  $v$ is a map $\AtProp \longrightarrow \P X$.
\end{definition}
The advantage of moving from the language of two-sorted frames to the language of enriched two-sorted frames is that the closure operator $\cl$ defined on direct presentations (see \eqref{eq : def closure on sets}) can be expressed in the modal language associated with enriched two-sorted frames. 
Indeed, unravelling the definitions involved, it is not difficult to see that for each 
subset $S$,
\begin{equation}
\label{eq : enriched two sorted frames and closure}
\cl(S) = \overline{{\downarrow}_{\leq_J}S} = \ef {\downarrow}_{\leq_J}S = \covD \niB \leqjD S.
\end{equation}
Recall that any assignment $v$ on a given finite lattice $L$ uniquely
gives rise to the  assignment $v^*$ on  $\mb{E}_L$ defined by
$v^*(p) := \{ j \in J(L) \mid j \leq_L v(p) \}$
for every $p \in AtProp$.  Then it can be readily verified that the following identity is satisfied for every $p\in \AtProp$:
$$v^*(p) = \covD \niB \leqjD v^*(p). $$
The semantic  identity above suggests the following definition:

\begin{definition}
A valuation $v$ on an  enriched two-sorted model $\mb{E}$ is \textit{closed} if 
$v(p) = \Drxy \Bryx \Drxx v(p)$ for every $p\in \AtProp$.
An  enriched two-sorted model is \textit{closed} if 
its associated valuation is closed.
\end{definition}

Thus, in the case a given enriched two-sorted model $\mb{M} = (\mb{E}_L,v)$  for some finite lattice $L$, 
the fact that the valuation $v$ arises from a lattice  assignment on $L$ can be expressed in the modal language of enriched two-sorted frames  by means of the satisfaction of the identity $p = \Drxy \Bryx \Drxx p$ for every $p \in \AtProp$.

\begin{definition}
An enriched two-sorted model is \textit{ordered} if its underlying  enriched two-sorted frame is ordered and its associated valuation assigns every $p\in \AtProp$ to a downset of $(X,\rxx)$.
\end{definition}

\subsection{Correspondence rules for enriched two-sorted frames}

In the present subsection, we show the soundness of the following extra rules on enriched two-sorted frames.

$$\frac{\Drxx \mbf{j} \W s
\leq \kappa(\mbf{j})}{s
\leq \kappa(\mbf{j})} \leqAtom$$
\begin{lemma}
\label{lem : soundness leq atom}
The rule $\leqAtom$ is sound and invertible on ordered enriched two-sorted models.
\end{lemma}
\begin{proof}
Fix an ordered enriched two-sorted model $\mb{M} = (\mb{E},v)$.
Let $x\in X$ such that   $v(\nomj)= \{ x\}$, and assume that 
$\Drxx \mbf{j} \W s
\leq \kappa(\mbf{j})$ is satisfied on $\mb{M}$. 
This means that $ x\notin ( {\downarrow_{\rxx}} x) \cap v(s) $, which implies that
$x\notin v(s)$. This condition is equivalent to $s
\leq \kappa(\mbf{j})$ being satisfied.

Conversely, assume that $s \leq \kappa(\mbf{j})$ is satisfied on $\mb{M}$ as above. 
This is equivalent to $x \notin v(s)$. Since by assumption $v(s)$ is a downset of $(X , \rxx)$, we have that $x \notin v(s)$ iff
$ {\downarrow_{\rxx}} x \nsubseteq v(s)$. Hence $x \notin {\downarrow_{\rxx}} x \cap v(s)$, which is equivalent to 
$\Drxx \mbf{j} \W s
\leq \kappa(\mbf{j})$ being satisfied on $\mb{M}$, as required.
\end{proof}
$$
\frac{\nomj \leq \covD \mbf{C} \quad \nomk \leq \inD \mbf{C}}{
\nomj \leq \covD \mbf{C} \quad \nomk \leq \inD \mbf{C}
\quad \covD \niB \leqxD( \inD\mbf{C} \smallsetminus \nomk ) \leq \kappa(\nomk)}
\MinCovB
$$
\begin{lemma}
\label{lem : MinCovB}
The rule $\MinCovB$ is sound and invertible on
every closed model $\mb{M} = (\mb{E}_L, v)$ such that $\mb{E}_L = ( J(L), \P J(L), \cov,\ni, \leq_J  )$ is the enriched two-sorted frame 
associated with
some finite lattice $L$ 
(cf.\ Definition \ref{def : two sorted frame of a lattice}).
\end{lemma}
\begin{proof}
The direction from bottom to top is immediate. Conversely, 
assume that the inequalities 
$\nomj \leq \covD \mbf{C} $ and $ \nomk \leq \inD \mbf{C}$
are satisfied on $\mb{M}$.
Let $j,k\in J(L)$ and $C\subseteq J(L)$ such that 
$v(\nomj)=\{j\}$, $v(\nomk)=\{k\}$ and $v(\mbf{C})=\{C\}$.
Hence,
$C\in \M(j)$ and $k\in C$. By Lemma \ref{lem : minimal cover closure}.2, this implies  that 
$k \notin \overline{{\downarrow_{\leq}}( C \smallsetminus k)}$, 
which is equivalent to 
the satisfaction of the inequality 
$\covD \niB \leqxD( \inD\mbf{C} \smallsetminus \nomk ) \leq \kappa(\nomk)$ on $\mb{M}$.
\end{proof}
\begin{lemma}
\label{lem : MinCovD}
Let $s$ be a $\mc{L}^+$-term.
For every closed model $\mb{M} = (\mb{E}_L, v)$ such that $\mb{E}_L = ( J(L), \P J(L), \cov,\ni, \leq_J  )$ is the enriched two-sorted frame 
associated with
some finite lattice $L$ 
(cf.\ Definition \ref{def : two sorted frame of a lattice}),
$$\mb{M} \Vdash (S1) \quad \tiff \quad \mb{M} \Vdash (S2), $$
where
\begin{align*}
    & (S1) := \sys{%
	    & \nomj \leq \covD \mbf{C} \\
        &  \nomk \leq \inD \mbf{C}\\
        & \leqjD \nomj \W \leqjD \nomk \leq \kappa(\nomk)\\
        & \leqjD \nomk \W s \leq \kappa(\nomk)\\
    },
\\
&   (S2):= \sys{%
	    & \nomj \leq \covD \mbf{C} \\
        &  \nomk \leq \inD \mbf{C}\\
        & \leqjD \nomj \W \leqjD \nomk \leq \kappa(\nomk)\\
        & \leqjD \nomk \W s \leq \kappa(\nomk)\\
        & \nomj \W \covD \niB   ( \covD\niB\leqjD ( \inD \mbf{C} \smallsetminus \nomk  ) \V   (\leqjD \nomj \W \leqjD \nomk)  
\V   (\leqjD \nomk \W s)) \leq \bot\\
    }.
\end{align*}
\end{lemma}
\begin{proof}
The right to left direction is immediate, since $(S1)$ is a  subset of $(S2)$. Assume that $(S1)$ is satisfied on $\mb{M}$. 
Let $j,k\in J(L)$ and $C\subseteq J(L)$ such that 
$v(\nomj)=\{j\}$, $v(\nomk)=\{k\}$ and $v(\mbf{C})=\{C\}$.
The assumptions imply that
$$C\in \M(j), \quad k\in C, \quad  
k \notin {\downarrow_{\leq_J}} j \cap {\downarrow_{\leq_J}} k,
\quad
k \notin {\downarrow_{\leq_J}} k \cap v(\phi).
$$
It is enough to show that 
$$ j \notin v(\covD \niB   ( \covD\niB\leqjD ( \inD \mbf{C} \smallsetminus \nomk  ) \V   (\leqjD \nomj \W \leqjD \nomk)  
\V   (\leqjD \nomk \W s))). $$
Unravelling the definitions of $\covD$ and $\niB$, the condition above is equivalent  to the following:
\begin{equation}
\label{eq : condition closure j proof mincovD}
\text{there exists no } D \subseteq J(L) \text{ such that } D\in \M(j) \tand
D \subseteq   v'(\covD\niB\leqjD ( \inD \mbf{C}  \smallsetminus  \nomk )  \V
(\leqjD \nomj\W \leqjD \nomk) \V (\leqjD \nomk \W s)). 
\end{equation}
The conditions $k \notin {\downarrow_{\leq_J}} j \cap {\downarrow_{\leq_J}} k$ and 
$k \notin {\downarrow_{\leq_J}} k \cap v(\phi)$
respectively imply that 
$k \notin {\downarrow_{\leq_J}} j $ and $k\notin v'(s)$.
Hence, the following chain of inclusions holds:
\begin{align*}
 & \; v'(\covD\niB\leqjD ( \inD \mbf{C}  \smallsetminus  \nomk )  \V
(\leqjD \nomj\W \leqjD \nomk) \V (\leqjD \nomk \W s))
\\
 = & \; v'( \covD\niB\leqjD ( \inD \mbf{C}  \smallsetminus  \nomk ) ) \cup ( v'
(\leqjD \nomj) \cap v'( \leqjD \nomk)) \cup (v' (\leqjD \nomk )\cap v'( s))
\\
 = & \; v'( \covD\niB\leqjD ( \inD \mbf{C}  \smallsetminus  \nomk ) ) \cup ( {\downarrow_{\leq_J}} j \cap {\downarrow_{\leq_J}} k) \cup ({\downarrow_{\leq_J}} k \cap v'( s))
\tag{by definition of $\leqjD$}
\\
\subseteq   & \; v'( \covD\niB\leqjD ( \inD \mbf{C}  \smallsetminus  \nomk ) ) \cup ( {\downarrow_{\leq_J}}k \smallsetminus \{k\})  \cup ({\downarrow_{\leq_J}} k \cap v'( s))
\tag{ $ k\notin {\downarrow_{\leq_J}} j $}
\\
\subseteq   & \; v'( \covD\niB\leqjD ( \inD \mbf{C}  \smallsetminus  \nomk ) ) \cup ( {\downarrow_{\leq_J}}k \smallsetminus \{k\})  \cup ({\downarrow_{\leq_J}} k \smallsetminus \{k\})
\tag{ $ k\notin v'(s)$}
\\
\subseteq   & \; v'( \covD\niB\leqjD ( \inD \mbf{C}  \smallsetminus  \nomk ) ) \cup ( {\downarrow_{\leq_J}}k \smallsetminus \{k\}) 
\\
\subseteq   & \; \overline{{\downarrow_{\leq_J}}(C \smallsetminus k)} \cup ( {\downarrow_{\leq_J}}k \smallsetminus \{k\}) .
\tag{cf.\ \eqref{eq : enriched two sorted frames and closure}}
\end{align*}
The fact that $C\in \M(j)$ and $k \in C$ implies, by
Lemma \ref{lem : minimal cover closure}.4, that
there is no cover of $j$ which is included in the set $\overline{{\downarrow_{\leq_J}}(C \smallsetminus k)} \cup  ( {\downarrow_{\leq_J}}k \smallsetminus \{k\}) $. 
This implies that \eqref{eq : condition closure j proof mincovD} holds as required.
\end{proof}

\subsection{Closed right Ackermann rule}
In the present section, we are going to prove the soundness of the following special version of the Ackermann rule.

$$\frac{ \exists p \; \Big [ 
\mbox{\Large{\&}}_{i=1}^{n} \{ \alpha_i \leq p   \} \;\; \& \;\;
 \mbox{\Large{\&}}_{j=1}^m \{ \beta_j(p)\leq \gamma_j(p)\} \Big ]
}{ 
\mbox{\Large{\&}}_{j=1}^m 
\{\beta_j(\covD  \niB \leqxD \bigvee_{i=1}^n \alpha_i)\leq \gamma_j(\covD  \niB \leqxD \bigvee_{i=1}^n \alpha_i) \} } (\ackRcl)$$

\begin{lemma}[Right Ackermann Lemma for closed models]
The rule $(\ackRcl)$ is sound and invertible on closed  enriched two-sorted models $\mb{M} = (\mb{E}_L, v)$ such that $\mb{E}_L$ is the enriched two-sorted frame associated with some finite lattice $L$.
%
\end{lemma}

\begin{proof}
Fix a closed  enriched two-sorted model $\mb{M} = (\mb{E}_L, v)$ such that $\mb{E}_L$ is the enriched two-sorted frame associated with some finite lattice $L$.
For the direction from bottom to top, 
assume that for every $1\leq j\leq m$,
$$ v (
\beta_j(\covD  \niB \leqxD \bigvee_{i=1}^n \alpha_i) )
\subseteq 
v( \gamma_j(\covD  \niB \leqxD \bigvee_{i=1}^n \alpha_i)). $$
Let $v'$ be  the $p$-variant of $v$ such that 
$v'(p) = \covD  \niB \leqxD v( \bigvee_{i=1}^n \alpha_i )$.
As discussed in subsection \ref{subsec : Enriched two sorted frame}, the composition $\covD \niB \leqxD$ is the operator that maps each set to the closure of its downset. Hence $v'(p)$ is a closed set and $v'$ is a closed valuation.
Since $\alpha_i$ does not contain $p$, we have that $v'(\alpha_i)=v(\alpha_i)$, and hence  
$$v'(\alpha_i) \leq v(\alpha_1) \V ... \V v(\alpha_n) \leq \covD  \niB \leqxD (v(\alpha_1) \V ... \V v(\alpha_n)) = v'(p).$$ This shows
that $$v'(\alpha_i) \leq v'(p), \text{ for all } 1\leq i \leq n.$$ 
Moreover, for all $1\leq i\leq m$, we have $$v'(\beta_i(p)) = v(\beta_i(\covD  \niB \leqxD (\alpha_1 \V ... \V \alpha_n )/p)) \leq v(\gamma_i(\covD  \niB \leqxD  (\alpha_1 \V ... \V \alpha_n) / p )) = v'(\gamma_i( p)).$$

For the implication from top to bottom, we make use of the fact that the $\beta_i$ are monotone (since positive) in $p$, while the $\gamma_i$ are antitone (since negative) in $p$.
Since the $\alpha_i$ do not contain $p$, and $v$ is a $p$-variant of $v'$, we have $v(\alpha_i)= v'(\alpha_i)\leq v'(p)$, for all $1\leq i \leq n$; hence, $v(\alpha_1) \V  ... \V v(\alpha_n) \leq v'(p)$. Since $v'$ is a closed valuation, $v'(p)$ is a closed set, and 
$v(\alpha_1) \V  ... \V v(\alpha_n) \leq v'(p)$ implies that 
$$\covD  \niB \leqxD (v(\alpha_1) \V  ... \V v(\alpha_n)) \leq v'(p).$$
Hence, 
$$ v(\beta_i(\covD  \niB \leqxD (\alpha_1 \V ... \V \alpha_n) / p ) )\leq
v'(\beta_i(p))\leq v'(\gamma_i(p)) \leq v(\gamma_i(\covD  \niB \leqxD  (\alpha_1 \V ... \V \alpha_n) / p ) ) .$$
\end{proof}

\section{Characterizing uniform upper bounds on the length of $\D$-chains in finite lattices }
\label{sec : Dchain}

\begin{definition}
\label{def : Dchain}
    Let $L$ be a finite lattice and let $\mb{E}_L = (J(L), \P J(L), \cov, \ni, \leq_J )$ be its associated enriched two-sorted frame (cf.\ Definition \ref{def : two sorted frame of a lattice}).
    Consider the  binary relation $\D \subseteq J(L) \times J(L)$ defined as follows: for any $j,k \in J(L)$, 
    $$ j \D k \quad\tiff\quad j \cov C, \;  k\in C \tand k \nleq j \text{ for some }  C \in \P J(L).$$
    A \textit{$\D$-chain} of length $l$ is a sequence $(j_0, ..., j_l)$ of elements of $J(L)$ 
    such that $j_i \D j_{i+1}$ for each $0 \leq i \leq (l-1)$.
\end{definition}


A notion similar to the one defined above has been used
in \cite{NationApproach} and \cite{SemenovaSuborders} to define a hierarchy of varieties of lattices progressively generalising the variety of distributive lattices. More discussion about the similarities and differences between Nation and Semenova's notion of $D$-chain and the one above can be found in Section \ref{sec : Conclusions}. In \cite{duality}, the result in \cite{NationApproach} and \cite{SemenovaSuborders}
has been generalized, and the existence of a Sahlqvist-type correspondence mechanism underlying it has been observed. The main motivation of the present paper is to provide a formal framework where this observation can be precisely spelled out, and 
 in the present section, we are ready to obtain  a result similar to Nation's  by means of an ALBA  reduction.

\bigskip

Fix enumerations of variables $x_n,y_n $ for $n\in \mb{N}$.
Consider
the following family of lattice inequalities:
$$\{ t_n \leq s_n  \mid n\in \mb{N}\},$$
such that the lattice terms $t_n$ and $s_n$ 
are recursively defined as follows:

\begin{center}
\begin{tabular}{cl}
  $t_0 \: := \: x_0$
  & $t_{n+1} \: := \: x_{n+1} \W 
  ( y_{n+1} \V t_n)$ \\
  $s_0 \: := \: \bot$
  & $s_{n+1} \: := \:  \: x_{n+1} \W 
  ( y_{n+1} \V  (x_{n+1} \W x_n ) \V s_n).$
  \\ & 
\end{tabular}
\end{center}

The aim of this section is proving the following proposition:
\begin{proposition}
\label{prop : Dchain}
For any finite lattice $L$ and any $n\in \mb{N}$, 
$$ L \models t_n \leq s_n \quad \tiff \quad  \text{ there is  no }
\D\text{-chain of length } n \text{ in } L. $$
\end{proposition}

\begin{proof}

For $n=0$, we need to prove that, if $L$ is a finite lattice,
$L \models x_0 \leq \bot$ iff   there is  no $\D$-chain of length $0$ in $L$.
This is clear, since the only finite lattice $L$  
such that $L \models x_0 \leq \bot$
is the one-element lattice, which is 
the only finite lattice which has no join-irreducible element.

\medskip

Let $n+1 \geq 1$. By Corollary \ref{cor : lattice frame equi}, 
$$L \models t_{n+1} \leq s_{n+1}
\quad \tiff \quad 
 \mb{E}_L \Vdash ST(t_{n+1}) \leq ST(s_{n+1}) ,$$
where $\mb{E}_L = (J(L), \P J(L) , \cov, \ni, \leq_J)$
is the enriched two-sorted frame associated with $L$ (cf.\ Definition \ref{def : two sorted frame of a lattice}), and
the validity on the right-hand side of the equivalence is understood in terms of satisfaction for every closed valuation.


We will
provide  ALBA$^{l}$ reductions
for each $n+1\geq 1$ and each inequality $t_{n+1} \leq s_{n+1}$.
Since all the ALBA$^{l}$ rules are sound and invertible on $\mb{E}_L$, the reduction will output a condition in the first order correspondence language of $\mb{E}_L$, 
which is equivalent to the validity of the input inequality on $L$, and 
which will express the existence of no $D$-chains of length $n+1$ in $L$.

First of all, using the standard translation introduced 
in section \ref{subsec : MML},
the lattice terms $t_{n+1}$ and $s_{n+1}$ translate into the following monotone modal logic formulas:
\begin{center}
\begin{tabular}{cll}
    $ST(t_0) \: := \: x_0$
    & $ST(t_{n+1}) \: := \: x_{n+1} \W 
    t_{n+1}'$ &  with  
    $t_{n+1}' := \ef ( y_{n+1} \V ST(t_n))$ \\
     $ST(s_0) \: := \: \bot$
    & $ST(s_{n+1}) \: := \:  \: x_{n+1} \W 
    s_{n+1}'$
    & with $s_{n+1}' := \ef ( y_{n+1} \V  (x_{n+1} \W x_n ) \V ST(s_n)).$
\end{tabular}
\end{center}

Using the notation introduced in section \ref{subsec : MML},
the $\mc{L}_{MML}$-terms above can be translated 
into the modal language of enriched two-sorted 
frames (cf.\ Subsection \ref{subsec : Enriched two sorted frame}) as indicated below.
For the sake of simplicity, we use the symbols $t_n$ and $s_n$ also to indicate the translations of the original lattice terms.
\label{def sn tn}


\begin{center}
\begin{tabular}{ccll}
$t_0 \: = \: x_0$ 
& $\quad$
& $t_{n+1} \: = \: x_{n+1} \W t_{n+1}'$
& with $t_{n+1}' \: = \:\covD  \niB ( y_{n+1} \V t_n) $\\
$s_0 \: = \: \bot$
& $\quad$
& $s_{n+1} \: = \:  \: x_{n+1} \W s_{n+1}'$
& with $s_{n+1}' \: = \: 
\covD  \niB ( y_{n+1} \V  (x_{n+1} \W x_n ) \V s_n) $
\end{tabular}
\end{center}

Let
$\overline{x}$ stand for the list of variables $x_{n} ,... ,x_0$, and
$\overline{y}$ stand for the list of variables $y_{n} ,..., y_1$.
ALBA$^l$ transforms the input inequality $t_{n+1} \leq s_{n+1}$
into the following quasi-inequality (cf.\ section \ref{subsec : intro ALBA}):
\begin{align*}
    \forall x_{n+1},\forall \overline{x}, 
    \forall y_{n+1},\forall \overline{y},
    \forall \mbf{j_{n+1}},\;
    & \sys{\sys{%
    & \mbf{j_{n+1}} \leq t_{n+1} \\
    & s_{n+1}  \leq \kappa(\mbf{j_{n+1}})\\
  }%
  \implies \false }.
\end{align*}
Since $t_{n+1} = x_{n+1} \W t_{n+1}'$ and $s_{n+1} = x_{n+1} \W s_{n+1}'$, we can rewrite the quasi-inequality above as:
\begin{align*}
    \forall x_{n+1},\forall \overline{x}, 
    \forall y_{n+1},\forall \overline{y},
    \forall \mbf{j_{n+1}},\;
    & \sys{\sys{%
	    & \mbf{j_{n+1}} \leq x_{n+1} \W 
	    t_{n+1}' \\
	    & x_{n+1} \W 
	    s_{n+1}'
	     \leq \kappa(\mbf{j_{n+1}})\\
    }
  \implies \false }.
\end{align*}
%
Applying the rule $(\spW )$ to the first inequality yields:
\begin{align*}
    \forall x_{n+1},\forall \overline{x}, 
    \forall y_{n+1},\forall \overline{y},
    \forall \mbf{j_{n+1}},\;
    & \sys{\sys{%
	    & \mbf{j_{n+1}} \leq x_{n+1} 
	    \\ &
	     \mbf{j_{n+1}} \leq t_{n+1}' 
	    \\ & 
	     x_{n+1} \W 
	    s_{n+1}'
	     \leq \kappa(\mbf{j_{n+1}})\\
    }
  \implies \false }.
\end{align*}
Notice that $x_{n+1} \notin Var(t_{n+1}')$ and $s_{n+1}'$ is monotone in $x_{n+1}$.
Thus we can apply the Ackermann rule 
 $(\ackRcl )$  to eliminate $x_{n+1}$ via the substitution
$x_{n+1} \longleftarrow \leqjD \mbf{j_{n+1}}$.
\begin{align*}
    \forall \overline{x}, 
    \forall y_{n+1},\forall \overline{y},
    \forall \mbf{j_{n+1}},\;
    & \sys{\sys{%
	    & \mbf{j_{n+1}} \leq t_{n+1}' \\
	    & \leqjD\mbf{j_{n+1}} \W 
	    s_{n+1}'(\leqjD\mbf{j_{n+1}}/x_{n+1})
	     \leq \kappa(\mbf{j_{n+1}})\\
    }%
  \implies \false }.
\end{align*}
Recall that $\mb{E}_L$ is an ordered enriched two-sorted frame and 
closed valuations assign variables to downsets.
Hence, by Lemma \ref{lem : soundness leq atom}, 
the quasi-inequality above is equivalent to the quasi-inequality below by applying the rule $(\leqAtom)$.
\begin{align*}
    \forall \overline{x}, 
    \forall y_{n+1},\forall \overline{y},
    \forall \mbf{j_{n+1}},\;
    & \sys{\sys{%
	    & \mbf{j_{n+1}} \leq t_{n+1}' \\
	    & s_{n+1}'(\leqjD\mbf{j_{n+1}}/x_{n+1})
	     \leq \kappa(\mbf{j_{n+1}})\\
    }%
  \implies \false }.
\end{align*}
By lemma \ref{lem : Induction},
the quasi-inequality above is equivalent to 
\begin{align}
\label{quasiInequality1}
    \forall \mbf{j_{n+1}},..., \mbf{j_0},
    \forall \mbf{C_n},\mbf{C_{n-1}},...,\mbf{C_0}\;
    & \sys{\sys{%
	    & \mbf{j_{n+1}} \leq \covD \mbf{C_{n}} \\
	    & \mbf{j_{n}} \leq \inD \mbf{C_{n}}\\
	    & \leqjD \mbf{j_{n+1}} \W \mbf{j_{n}} \leq \bot\\
	    & \ldots\\
	    & \mbf{j_1} \leq \covD \mbf{C_{0}} \\
	    & \mbf{j_{0}} \leq \inD \mbf{C_{0}}\\
	    & \leqjD \mbf{j_1} \W \mbf{j_{0}} \leq \bot \\
    }%
  \implies \false }.
\end{align}
Notice that, for $0\leq i\leq n$,
the following inequalities:
$$\mbf{j_{i+1}} \leq \covD \mbf{C_i}, \quad \mbf{j_i} \leq \inD \mbf{C_i}, \quad
\leqjD \mbf{j_{i+1}} \W \mbf{j_i} \leq \bot $$
are respectively equivalent  to the following atomic formulas
in the first order correspondence language of enriched two-sorted frames  (cf.\ Subsection \ref{subsec : Enriched two sorted frame}):
$$j_{i+1} \cov C_i,  \quad j_i \in C_i, \quad j_{i} \nleq j_{i+1}.$$
By definition \ref{def : Dchain}, the conditions above yield  $j_{i+1}\D j_i$ for each $0\leq i \leq n$. Hence the quasi-inequality (\ref{quasiInequality1})
is equivalent to the following quasi-inequality:
$$ \forall j_{n+1},... ,j_0\;\quad
    \lb 
    (\;  j_{n+1} \D j_n \ldots j_1 \D j_{0} \; )
	\implies \false \rb,$$
which expresses the condition that
there is no $\D$-chain of length $ n+1$.
\end{proof}

The proof of the proposition above relies on the following lemma,
the proof of which can be found in section \ref{subsec : proof of lemma induction}.

\begin{lemma}
\label{lem : Induction}
    For every $n\geq 1$, 
      ALBA$^l$ succeeds
    on the quasi-inequality
    \begin{align*}
    \forall x_{n-1},..., x_0,
    \forall y_n,...,y_0,
    \forall \mbf{j_n},\;
    & \sys{\sys{%
	    & \mbf{j_n} \leq t_n' \\
	    & s_n'( \leqjD\mbf{j_n} / x_n )\leq \kappa(\mbf{j_n}) \\
    }%
  \implies \false },
\end{align*}
and produces
\begin{align*}
    \forall \mbf{j_n},... \mbf{j_0},
    \forall \mbf{C_{n-1}},...\mbf{C_0}\;
    & \sys{\sys{%
	    & \mbf{j_n} \leq \covD \mbf{C_{n-1}} \\
	    & \mbf{j_{n-1}} \leq \inD \mbf{C_{n-1}}\\
	    & \leqjD \mbf{j_n} \W \mbf{j_{n-1}} \leq \bot\\
	    & \ldots\\
	    & \mbf{j_1} \leq \covD \mbf{C_{0}} \\
	    & \mbf{j_{0}} \leq \inD \mbf{C_{0}}\\
	    &\leqjD \mbf{ j_1} \W \mbf{j_{0}} \leq \bot \\
    }%
  \implies \false }.
\end{align*}
\end{lemma}


\section{Conclusions and further directions }
\label{sec : Conclusions}
    \paragraph{Conclusions.} In the present paper, the algorithmic correspondence theory revolving around ALBA (cf.\ \cite{ALBA, Gilardi}) has been adapted and extended, so as to provide an adequate environment in which to formalize the observation (cf.\ \cite{duality}) of the existence of a Sahlqvist-type mechanism underlying 
dual characterization results
for finite lattices.

The treatment of lattice inequalities in the setting of ALBA is mediated by monotone modal logic, thanks to the existence of a duality-on-objects between finite lattices and  join-presentations (cf.\ Definition \ref{def : join presentation}), and the fact that join-presentations are closely related to (monotone) neighbourhood frames.

A key step towards the main result of the present paper is the adaptation of ALBA to monotone modal logic, semantically justified by the introduction of two-sorted structures and their associated correspondence language. In this setting, the Sahlqvist correspondence theory of \cite{HelleThesis} can be embedded and generalized.

\paragraph*{Comparison with Nation's results.}
As mentioned early on, 
our result is similar to Nation's dual characterization  of uniform upper bounds on the length of $D$-chains in finite lattices:

\begin{definition}
    Let $L$ be a finite lattice.
    Let $D \subseteq J(L) \times J(L)$ be the  binary relation defined as follows:  for any $j,k \in J(L)$, 
    $$ j D k \quad\tiff\quad j \cov C, \;  k\in C \tand k \neq j \text{ for some }  C \in \P J(L).$$
    A \textit{$D$-chain} of length $l$ is a sequence $(j_0, ..., j_l)$ of elements of $J(L)$ 
    such that $j_i D j_{i+1}$ for each $0 \leq i \leq (l-1)$.
\end{definition}

The dual characterization of section \ref{sec : Dchain} is different  
from Nation's \cite{NationApproach} and is not covered by the result in \cite{duality}, which generalizes Nation's. As far as we know, it is original.

Clearly, $\D$ is included in $D$ for any finite lattice $L$.
Hence, the validity on a finite lattice $L$ of 
Nation's inequalities for a given $n$ is a sufficient condition  for $L$ having $\D$-chains of length at most $n$.
However, in the remainder of the paragraph,
we are going to show
that 
this upper bound is not accurate.
Indeed,
the maximal length of $\D$-chains starting from a given join-irreducible element in a lattice  can be strictly smaller than the one for $D$-chains starting from the same join-irreducible element.
Consider the lattice $L$ the Hasse diagram of which is given by the figure. In this example, it can be easily verified that
\begin{align*}
J(L) & := \set{ a,b,c,d,e},\\
\M(c)& := \set{\set{c}, \set{a,b} },\\
\M(e) & := \set{ \set{e}, \set{a,d},\set{c,d} }.
\end{align*}
The only  $\D$-chain starting from $e$ is $e\D d$, whereas there are $D$-chains of length 2 starting  from $e$, for instance $eDcDb$.

\begin{figure}[h]
  \begin{center}
    \includegraphics[scale=0.4]{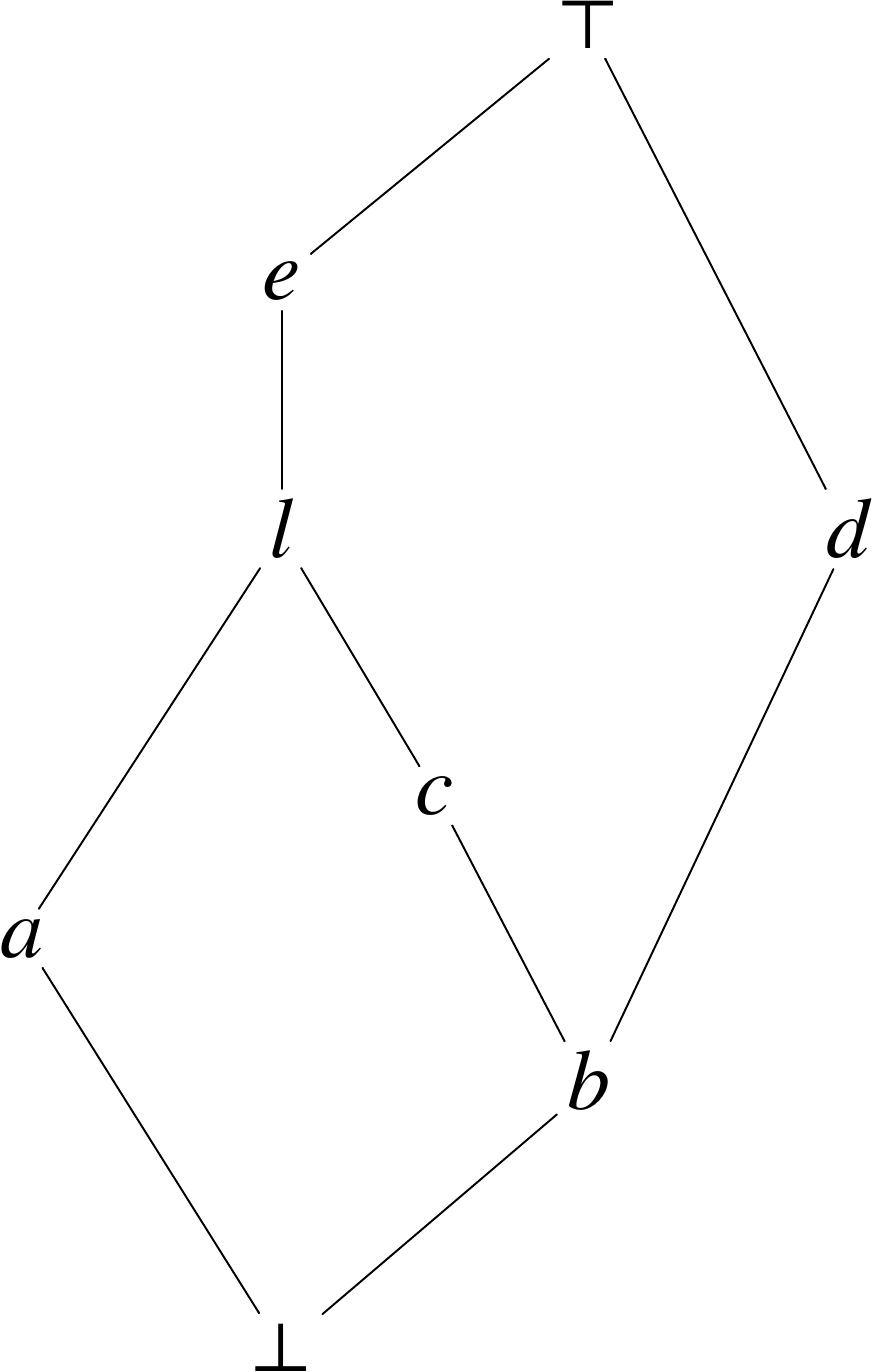}
  \end{center}
  \caption{The Hasse diagram of the lattice $L$}
  \label{Fig : lattice}
\end{figure}

\paragraph{Further directions.} 
The present paper is a first step towards the fully-fledged automatization  of dual characterization results for finite lattices. Significant extensions of Nation's dual characterization results appear e.g.\ in \cite{SemenovaSuborders} and \cite[Proposition 8.5]{duality}. Hence, natural directions worth pursuing are (a) extending the results of the present paper so as to account for  \cite[Proposition 8.5]{duality}, and (b)  analyzing  the technical machinery introduced in the present paper from an algorithmic perspective. The latter point involves e.g.\ establishing whether the present set of rules is minimal, or whether some rules can actually be derived.

Related to both these directions, but more on the front of methodology, are outstanding open questions about Lemma \ref{lem : double ackermann}. This lemma provides the soundness and invertibility of a rule by means of which variable elimination is effected via instantiation. So far, all rules of this type in ALBA have been proved sound and invertible thanks to one or another  version of Ackermann's lemma. However, it is not clear whether Lemma \ref{lem : double ackermann} can be  accounted for in terms of Ackermann's lemma, and hence whether the rule justified in Lemma \ref{lem : double ackermann} can be regarded as an Ackermann-type rule. Moreover, while Lemma \ref{lem : double ackermann} is rooted and has an intuitive understanding in the semantics of minimal coverings, at the moment it is not clear whether and how more general versions of  this rule  can be formulated, which would be of a wider applicability. Giving answers to these questions would significantly enlarge the scope of  algorithmic correspondence theory, and is also a worthwhile future direction.

\bibliographystyle{plain}
\bibliography{2sorted}

\newpage

\appendix
\label{Appendix}

\section{Proof of lemma \ref{lem : Induction}}
\label{subsec : proof of lemma induction}
    
\textbf{Lemma \ref{lem : Induction}}
    For every $n\geq 1$, 
      ALBA$^l$ succeeds
    on the quasi-inequality
    \begin{align}
    \label{initialSystem}
    \forall x_{n-1},..., x_0,
    \forall y_n,...,y_0,
    \forall \mbf{j_n},\;
    & \sys{\sys{%
	    & \mbf{\mbf{j_n}} \leq t_n' \\
	    & s_n'( \leqxD\mbf{\mbf{j_n}} / x_n )\leq \kappa(\mbf{\mbf{j_n}}) \\
    }%
  \implies \false },
\end{align}
and produces
\begin{align*}
    \forall \mbf{\mbf{j_n}},... \mbf{\mbf{j_0}},
    \forall \mbf{C_{n-1}},...\mbf{C_0}\;
    & \sys{\sys{%
	    & \mbf{\mbf{j_n}} \leq \covD \mbf{C_{n-1}} \\
	    & \mbf{j_{n-1}} \leq \inD \mbf{C_{n-1}}\\
	    & \leqjD \mbf{\mbf{j_n}} \W \mbf{j_{n-1}} \leq \bot\\
	    & \ldots\\
	    & \mbf{j_1} \leq \covD \mbf{C_{0}} \\
	    & \mbf{j_{0}} \leq \inD \mbf{C_{0}}\\
	    & \leqjD \mbf{j_1} \W \mbf{j_{0}} \leq \bot \\
    }%
  \implies \false }.
\end{align*}
%
\begin{proof}
By induction on $n$. 
If $n=1$, then the quasi-inequality (\ref{initialSystem})
has the following shape (cf.\ the definitions of $t_n$ and $s_n$ on page \pageref{def sn tn}):
\begin{align*}
    \forall  x_0,
    \forall y_1,
    \forall \mbf{j_1},\;
    & \sys{\sys{%
	    & \mbf{j_1} \leq \covD  \niB ( y_{1} \V x_0) \\
	    & \covD  \niB ( y_{1} \V  (\leqjD\mbf{j_{1}} \W x_0 ) \V \bot) 
	    \leq \kappa(\mbf{j_1}) \\
    }%
  \implies \false } .
\end{align*}
By applying  the rules  $( \AtCoat )$ and $( \tWB )$ to the second inequality, we get:
\begin{align*}
    \forall  x_0,
    \forall y_1,
    \forall \mbf{j_1},\;
    & \sys{\sys{%
	    & \mbf{j_1} \leq \covD  \niB ( y_{1} \V x_0) \\
	    & \mbf{j_1} \leq \neg \covD  \niB ( y_{1} \V  (\leqjD\mbf{j_{1}} \W x_0 ) ) \\
    }%
  \implies \false } .
\end{align*}
Now we can apply $( \tDB )$ and $(\tBD )$ to the second inequality, and get:
\begin{align*}
    \forall  x_0,
    \forall y_1,
    \forall \mbf{j_1},\;
    & \sys{\sys{%
	    & \mbf{j_1} \leq \covD  \niB ( y_{1} \V x_0) \\
	    & \mbf{j_1} \leq  \covB  \niD \neg ( y_{1} \V  (\leqjD\mbf{j_{1}} \W x_0 ) ) \\
    }%
  \implies \false } .
\end{align*}
By applying  the rule  $(\tnm)$ to the first and second inequalities, we get:\begin{align*}
    \forall  x_0,
    \forall y_1,
    \forall \mbf{j_1},\;
    & \sys{\sys{%
	    & \mbf{j_1} \leq \covD  ( \niB ( y_{1} \V x_0) \W 
	    \niD \neg ( y_{1} \V  (\leqjD\mbf{j_{1}} \W x_0 ) ) ) \\
	    & \mbf{j_1} \leq  \covB  \niD \neg ( y_{1} \V  (\leqjD\mbf{j_{1}} \W x_0 ) ) \\
    }%
  \implies \false } .
\end{align*}
We can apply the rule  $(\apD)$ to the first inequality, and get:
\begin{align*}
    \forall  x_0,
    \forall y_1,
    \forall \mbf{j_1}, \forall \mbf{C_0} \;
    & \sys{\sys{%
	    & \mbf{j_1} \leq \covD \mbf{C_0}\\
	    &\mbf{C_0} \leq   \niB ( y_{1} \V x_0) \W 
	    \niD \neg ( y_{1} \V  (\leqjD\mbf{j_{1}} \W x_0 ) )  \\
	    & \mbf{j_1} \leq  \covB  \niD \neg ( y_{1} \V  (\leqjD\mbf{j_{1}} \W x_0 ) ) \\
    }%
  \implies \false } .
\end{align*}
By applying  the rule  $(\spW)$ to the second inequality, we get:
\begin{align*}
    \forall  x_0,
    \forall y_1,
    \forall \mbf{j_1}, \forall \mbf{C_0} \;
    & \sys{\sys{%
	    & \mbf{j_1} \leq \covD \mbf{C_0}\\
	    &\mbf{C_0} \leq   \niB ( y_{1} \V x_0) \\
	    &\mbf{C_0} \leq 
	    \niD \neg ( y_{1} \V  (\leqjD\mbf{j_{1}} \W x_0 ) )  \\
	    & \mbf{j_1} \leq  \covB  \niD \neg ( y_{1} \V  (\leqjD\mbf{j_{1}} \W x_0 ) ) \\
    }%
  \implies \false } .
\end{align*}
We can now apply  the rule  $(\ajB)$ to the second inequality and $(\apD)$ to the third inequality, and get:
\begin{align*}
    \forall  x_0,
    \forall y_1,
    \forall \mbf{j_1}, \mbf{j_0}, \forall \mbf{C_0} \;
    & \sys{\sys{%
	    & \mbf{j_1} \leq \covD \mbf{C_0}\\
	    &\inD \mbf{C_0} \leq     y_{1} \V x_0 \\
	    &\mbf{C_0} \leq 
	    \niD \mbf{j_0} \\
	    & \mbf{j_0} \leq 
	    \neg ( y_{1} \V  (\leqjD\mbf{j_{1}} \W x_0 ) )  \\
	    & \mbf{j_1} \leq  \covB  \niD \neg ( y_{1} \V  (\leqjD\mbf{j_{1}} \W x_0 ) ) \\
    }%
  \implies \false } .
\end{align*}
By applying  the rules  $( \DeMorgan )$ and $(\spW)$ to the fourth inequality, we get:
\begin{align*}
    \forall  x_0,
    \forall y_1,
    \forall \mbf{j_1}, \mbf{j_0}, \forall \mbf{C_0} \;
    & \sys{\sys{%
	    & \mbf{j_1} \leq \covD \mbf{C_0}\\
	    &\inD \mbf{C_0} \leq     y_{1} \V x_0 \\
	    &\mbf{C_0} \leq 
	    \niD \mbf{j_0} \\
	    & \mbf{j_0} \leq 
	    \neg  y_{1} \\
	    & \mbf{j_0} \leq  \neg  (\leqjD\mbf{j_{1}} \W x_0 ) )  \\
	    & \mbf{j_1} \leq  \covB  \niD \neg ( y_{1} \V  (\leqjD\mbf{j_{1}} \W x_0 ) ) \\
    }%
  \implies \false } .
\end{align*}
By applying 
the rule  $( \tRRminus )$  to the third inequality, 
the rules  $( \tWB )$ and $( \AtCoat )$  to the fourth and fifth inequalities, 
and the rules $( \tDB )$ and $(\tBD )$ to the last inequality, we get:
\begin{align*}
    \forall  x_0,
    \forall y_1,
    \forall \mbf{j_1}, \mbf{j_0}, \forall \mbf{C_0} \;
    & \sys{\sys{%
	    & \mbf{j_1} \leq \covD \mbf{C_0}\\
	    &\inD \mbf{C_0} \leq     y_{1} \V x_0 \\
	    &\mbf{j_0} \leq 
	    \inD \mbf{C_0} \\
	    & y_{1} \leq \kappa(\mbf{j_0} )\\
	    & \leqjD\mbf{j_{1}} \W x_0 \leq \kappa( \mbf{j_0} )\\
	    & \mbf{j_1} \leq \neg \covD  \niB  ( y_{1} \V  (\leqjD\mbf{j_{1}} \W x_0 ) ) \\
    }%
  \implies \false } .
\end{align*}
By applying  the rule  $( \tra )$  to the second and third inequalities and the rule $( \tWB )$  to the last inequality, we get:
\begin{align*}
    \forall  x_0,
    \forall y_1,
    \forall \mbf{j_1}, \mbf{j_0}, \forall \mbf{C_0} \;
    & \sys{\sys{%
	    & \mbf{j_1} \leq \covD \mbf{C_0}\\
	    &\inD \mbf{C_0} \leq     y_{1} \V x_0 \\
	    &\mbf{j_0} \leq 
	    \inD \mbf{C_0} \\
	    &\mbf{j_0} \leq     y_{1} \V x_0 \\
	    & y_{1} \leq \kappa(\mbf{j_0} )\\
	    & \leqjD\mbf{j_{1}} \W x_0 \leq \kappa( \mbf{j_0} )\\
	    &  \mbf{j_1} \W  \covD  \niB  ( y_{1} \V  (\leqjD\mbf{j_{1}} \W x_0 ) ) \leq \bot \\
    }%
  \implies \false } .
\end{align*}
%
%
By applying the rule $(\MT)$ to the fourth and fifth inequalities, 
the rule $(\AtCoat)$ to the last inequality,
and by exchanging the position of the second and third inequalities,
the quasi-inequality above can be equivalently rewritten as follows:
\begin{align*}
    \forall  x_0,
    \forall y_1,
    \forall \mbf{j_1},\mbf{j_0},\forall \mbf{C_0}\;
    & \sys{\sys{%
	    & \mbf{j_1} \leq \covD \mbf{C_0} \\
	    & \mbf{j_0} \leq \inD \mbf{C_0} \\
	    & \inD \mbf{C_0} \leq   y_{1} \V x_0 \\
	    & \mbf{j_0} \leq  x_0\\
	    &  y_{1} \leq \kappa(\mbf{j_0}) \\  
	    & \leqjD\mbf{j_{1}} \W x_0 \leq \kappa(\mbf{j_0}) \\
	    & \mbf{j_1} \W \covD  \niB ( y_{1} \V  (\leqjD\mbf{j_{1}} \W x_0 )) 
	    \leq \bot \\
    }%
  \implies \false } .
\end{align*}
By lemma \ref{lem : double ackermann} with the following instantiations\footnote{Notice that $t:=\top$ and $s:=\bot$ reduce the inequalities 
$\nomk \leq t$, $\inD \mbf{C}  \leq y \V t$ and 
$x \W s\leq \kappa(\nomk)$ in the statement of Lemma \ref{lem : double ackermann}  to tautologies, and the inequality
$\nomj \W 
	    \covD  \niB ( y \V  (\leqjD \nomj \W x ) \V (x \W s))\leq \bot$
to $\nomj \W \covD  \niB ( y \V  (\leqjD \nomj \W x ))\leq \bot$.}
$$t:=\top , \quad s:=\bot , \quad \nomj := \mbf{j_1} , \quad \nomk:= \mbf{j_0}
 , \quad \mbf{C}:=\mbf{C_0} , \quad x:=x_0,$$
the quasi-inequality above is equivalent to the following:
\begin{align*}
    \forall  x_0,
    \forall \mbf{j_1}, \mbf{j_0},
    \forall \mbf{C_0} \;
    & \sys{\sys{%
	    & \mbf{j_1} \leq \covD \mbf{C_0} \\
	    & \mbf{j_0} \leq \inD \mbf{C_0}\\
	    & \cl ( \inD \mbf{C_0}  \smallsetminus  \mbf{j_0} ) \leq \kappa(\mbf{j_0})\\
	    & \leqjD \mbf{j_1} \W \leqjD \mbf{j_0}  \leq \kappa(\mbf{j_0}) \\
	    & \mbf{j_1} \W 
	    \covD  \niB ( \cl ( \inD \mbf{C_0}  \smallsetminus  \mbf{j_0} ) 
	    \V (\leqjD \mbf{j_1} \W \leqjD \mbf{j_0} ))\leq \bot\\
    }%
  \implies \false } ,
\end{align*}
where $\cl$ abbreviates the composition $\covD \niB \leqjD$.
By applying the rule $\MinCovB$ bottom to top, the quasi-inequality above can be equivalently rewritten as follows
\begin{align*}
    \forall  x_0,
    \forall \mbf{j_1}, \mbf{j_0},
    \forall \mbf{C_0} \;
    & \sys{\sys{%
	    & \mbf{j_1} \leq \covD \mbf{C_0} \\
	    & \mbf{j_0} \leq \inD \mbf{C_0}\\
	    & \leqjD \mbf{j_1} \W \leqjD \mbf{j_0}  \leq \kappa(\mbf{j_0}) \\
	    & \mbf{j_1} \W 
	    \covD  \niB ( \cl ( \inD \mbf{C_0}  \smallsetminus  \mbf{j_0} ) 
	    \V (\leqjD \mbf{j_1} \W \leqjD \mbf{j_0} ))\leq \bot\\
    }%
  \implies \false } ,
\end{align*}
By applying Lemma \ref{lem : MinCovD}, we get:
\begin{align*}
    \forall  x_0,
    \forall \mbf{j_1}, \mbf{j_0},
    \forall \mbf{C_0} \;
    & \sys{\sys{%
	    & \mbf{j_1} \leq \covD \mbf{C_0} \\
	    & \mbf{j_0} \leq \inD \mbf{C_0}\\
	    & \leqjD \mbf{j_1} \W \leqjD \mbf{j_0}  \leq \kappa(\mbf{j_0}) \\
    }%
  \implies \false } ,
\end{align*}
By applying $(\leqAtom)$ to the third inequality, we get:
\begin{align*}
    \forall  x_0,
    \forall \mbf{j_1}, \mbf{j_0},
    \forall \mbf{C_0} \;
    & \sys{\sys{%
	    & \mbf{j_1} \leq \covD \mbf{C_0} \\
	    & \mbf{j_0} \leq \inD \mbf{C_0}\\
	    & \leqjD \mbf{j_1}  \leq \kappa(\mbf{j_0}) \\
    }%
  \implies \false } ,
\end{align*}
By $(\AtCoat)$ to the third inequality, we get:
\begin{align*}
    \forall  x_0,
    \forall \mbf{j_1}, \mbf{j_0},
    \forall \mbf{C_0} \;
    & \sys{\sys{%
	    & \mbf{j_1} \leq \covD \mbf{C_0} \\
	    & \mbf{j_0} \leq \inD \mbf{C_0}\\
	    & \leqjD \mbf{j_1} \W \mbf{j_0} \leq \bot \\
    }%
  \implies \false } ,
\end{align*}
which finishes the proof of the base case.

\paragraph*{Induction step.}
Fix $n\geq 1$, and assume that the lemma holds for $n$.
Recall that
$\overline{x}$ stands for the list of variables $x_{n} ,... ,x_0$, and
$\overline{y}$ stands for the list of variables $y_{n} ,..., y_1$.
Let us prove the lemma for the quasi-inequality
\begin{align*}
    \forall \overline{x},
    \forall y_{n+1},
    \forall \overline{y},
    \forall \mbf{j_{n+1}}\;
    & \sys{\sys{%
	    & \mbf{j_{n+1}} \leq t_{n+1}' \\
	    & s_{n+1}'( \leqjD \mbf{j_{n+1}} / x_{n+1} )\leq \kappa(\mbf{j_{n+1}}) \\
    }%
  \implies \false } .
\end{align*}
By the definitions on page \pageref{def sn tn},  the quasi-inequality above
can be rewritten into:
\begin{align*}
    \forall \overline{x},
    \forall y_{n+1},
    \forall \overline{y},
    \forall \mbf{j_{n+1}}\;
    & \sys{\sys{%
	    & \mbf{j_{n+1}} \leq \covD  \niB ( y_{n+1} \V t_n) \\
	    &\covD  \niB ( y_{n+1} \V  ( \leqjD \mbf{j_{n+1}} \W x_n ) \V s_n) \leq \kappa(\mbf{j_{n+1}})\\
    }
  \implies \false } ,
\end{align*}
which, by applying the rules $(\AtCoat )$ and $( \tWB )$ to the second inequality, is equivalent to:
\begin{align*}
    \forall \overline{x},
    \forall y_{n+1},
    \forall \overline{y},
    \forall \mbf{j_{n+1}}\;
    & \sys{\sys{%
	    & \mbf{j_{n+1}} \leq \covD  \niB ( y_{n+1} \V t_n) \\
	    & \mbf{j_{n+1}} \leq \neg 
	    \covD  \niB ( y_{n+1} \V  (\leqjD\mbf{j_{n+1}} \W x_n ) \V s_n) \\
    }
  \implies \false } .
\end{align*}
By applying  the rule  $( \tBD )$ and $(\tDB )$ to the second inequality, we get:
\begin{align*}
    \forall \overline{x},
    \forall y_{n+1},
    \forall \overline{y},
    \forall \mbf{j_{n+1}}\;
    & \sys{\sys{%
	    & \mbf{j_{n+1}} \leq \covD  \niB ( y_{n+1} \V t_n) \\
	    & \mbf{j_{n+1}} \leq  
	    \covB  \niD \neg ( y_{n+1} \V  (\leqjD\mbf{j_{n+1}} \W x_n ) \V s_n) \\
    }
  \implies \false } .
\end{align*}
By applying  the rule  $(\tnm)$ to the first and second inequalities, we get:
\begin{align*}
    \forall \overline{x},
    \forall y_{n+1},
    \forall \overline{y},
    \forall \mbf{j_{n+1}}\;
    & \sys{\sys{%
	    & \mbf{j_{n+1}} \leq \covD ( \niB ( y_{n+1} \V t_n)  \W 
	    \niD \neg ( y_{n+1} \V  (\leqjD\mbf{j_{n+1}} \W x_n ) \V s_n)  )\\
	    & \mbf{j_{n+1}} \leq  
	    \covB  \niD \neg ( y_{n+1} \V  (\leqjD\mbf{j_{n+1}} \W x_n ) \V s_n) \\
    }
  \implies \false } .
\end{align*}
By applying  the rule  $(\apD)$ to the first inequality, we get:
\begin{align*}
    \forall \overline{x},
    \forall y_{n+1},
    \forall \overline{y},
    \forall \mbf{j_{n+1}}, 
    \forall \mbf{C_{n}}\; 
    & \sys{\sys{%
	    & \mbf{j_{n+1}} \leq \covD \mbf{C_{n}}\\
	    & \mbf{C_{n}} \leq  \niB ( y_{n+1} \V t_n)  \W 
	    \niD \neg ( y_{n+1} \V  (\leqjD\mbf{j_{n+1}} \W x_n ) \V s_n)  \\
	    & \mbf{j_{n+1}} \leq  
	    \covB  \niD \neg ( y_{n+1} \V  (\leqjD\mbf{j_{n+1}} \W x_n ) \V s_n) \\
    }
  \implies \false } .
\end{align*}
By applying  the rule  $(\spW)$ to the second inequality, we get:
\begin{align*}
    \forall \overline{x},
    \forall y_{n+1},
    \forall \overline{y},
    \forall \mbf{j_{n+1}}, 
    \forall \mbf{C_{n}}\; 
    & \sys{\sys{%
	    & \mbf{j_{n+1}} \leq \covD \mbf{C_{n}}\\
	    & \mbf{C_{n}} \leq  \niB ( y_{n+1} \V t_n) \\
	    & \mbf{C_{n}} \leq
	    \niD \neg ( y_{n+1} \V  (\leqjD\mbf{j_{n+1}} \W x_n ) \V s_n)  \\
	    & \mbf{j_{n+1}} \leq  
	    \covB  \niD \neg ( y_{n+1} \V  (\leqjD\mbf{j_{n+1}} \W x_n ) \V s_n) \\
    }
  \implies \false } .
\end{align*}
We can now apply  the rule  $(\ajB)$ to the second inequality and $(\apD)$ to the third inequality, and get:
\begin{align*}
    \forall \overline{x},
    \forall y_{n+1},
    \forall \overline{y},
    \forall \mbf{j_{n+1}}, \mbf{j_{n}},
    \forall \mbf{C_{n}}\; 
    & \sys{\sys{%
	    & \mbf{j_{n+1}} \leq \covD \mbf{C_{n}}\\
	    & \inD \mbf{C_{n}} \leq   y_{n+1} \V t_n \\
	    & \mbf{C_{n}} \leq \niD \mbf{j_{n}} \\
	    &  \mbf{j_{n}} \leq \neg ( y_{n+1} \V  (\leqjD\mbf{j_{n+1}} \W x_n ) \V s_n)  \\
	    & \mbf{j_{n+1}} \leq  
	    \covB  \niD \neg ( y_{n+1} \V  (\leqjD\mbf{j_{n+1}} \W x_n ) \V s_n) \\
    }
  \implies \false } .
\end{align*}
By applying  the rules  $( \DeMorgan )$ and $(\spW)$ to the fourth inequality, we get:
\begin{align*}
    \forall \overline{x},
    \forall y_{n+1},
    \forall \overline{y},
    \forall \mbf{j_{n+1}}, \mbf{j_{n}},
    \forall \mbf{C_{n}}\; 
    & \sys{\sys{%
	    & \mbf{j_{n+1}} \leq \covD \mbf{C_{n}}\\
	    & \inD \mbf{C_{n}} \leq   y_{n+1} \V t_n \\
	    & \mbf{C_{n}} \leq \niD \mbf{j_{n}} \\
	    &  \mbf{j_{n}} \leq \neg  y_{n+1} \\
	    &  \mbf{j_{n}} \leq \neg  (\leqjD\mbf{j_{n+1}} \W x_n ) \\
	    &  \mbf{j_{n}} \leq \neg s_n  \\
	    & \mbf{j_{n+1}} \leq  
	    \covB  \niD \neg ( y_{n+1} \V  (\leqjD\mbf{j_{n+1}} \W x_n ) \V s_n) \\
    }
  \implies \false } .
\end{align*}
By applying  
the rule  $( \tRRminus )$  to the third inequality,
the rules  $( \tWB )$ and $( \AtCoat )$  to the fourth, fifth and sixth inequalities, 
and the rules $( \tDB )$ and $(\tBD )$ to the last inequality,
we get:
\begin{align*}
    \forall \overline{x},
    \forall y_{n+1},
    \forall \overline{y},
    \forall \mbf{j_{n+1}}, \mbf{j_{n}},
    \forall \mbf{C_{n}}\; 
    & \sys{\sys{%
	    & \mbf{j_{n+1}} \leq \covD \mbf{C_{n}}\\
	    & \inD \mbf{C_{n}} \leq   y_{n+1} \V t_n \\
	    & \mbf{j_{n}} \leq \inD \mbf{C_{n}} \\
	    &  y_{n+1} \leq \kappa(\mbf{j_{n}}) \\
	    &   \leqjD\mbf{j_{n+1}} \W x_n \leq \kappa(\mbf{j_{n}}) \\
	    &   s_n \leq \kappa(\mbf{j_{n}})  \\
	    & \mbf{j_{n+1}} \leq  \neg
	    \covD  \niB  ( y_{n+1} \V  (\leqjD\mbf{j_{n+1}} \W x_n ) \V s_n) \\
    }
  \implies \false } .
\end{align*}
By applying the rule $( \tra )$ to the second and third inequalities
  and the rule $( \tWB )$ to the last
inequality, and by exchanging the position of the second and third
inequalities, we get:
\begin{align*}
    \forall \overline{x},
    \forall y_{n+1},
    \forall \overline{y},
    \forall \mbf{j_{n+1}}, \mbf{j_{n}},
    \forall \mbf{C_{n}}\; 
    & \sys{\sys{%
	    & \mbf{j_{n+1}} \leq \covD \mbf{C_{n}}\\
	    & \mbf{j_{n}} \leq \inD \mbf{C_{n}} \\
	    & \inD \mbf{C_{n}} \leq   y_{n+1} \V t_n \\
	    & \mbf{j_{n}} \leq  y_{n+1} \V t_n\\
	    &  y_{n+1} \leq \kappa(\mbf{j_{n}}) \\
	    &   \leqjD\mbf{j_{n+1}} \W x_n \leq \kappa(\mbf{j_{n}}) \\
	    &   s_n \leq \kappa(\mbf{j_{n}})  \\
	    & \mbf{j_{n+1}} \W
	    \covD  \niB  ( y_{n+1} \V  (\leqjD\mbf{j_{n+1}} \W x_n ) \V s_n)\leq \bot \\
    }
  \implies \false } .
\end{align*}
By applying the rule $(\MT)$ to the fourth and fifth inequalities,
and since by definition $t_n = x_n \W t_n'$ and $s_n = x_n \W s_n'$, the quasi-inequality above is equivalent to the quasi-inequality below:
\begin{align*}
    \forall \overline{x},
    \forall y_{n+1},
    \forall \overline{y},
    \forall \mbf{j_{n+1}},
    \mbf{j_n},
    \forall \mbf{C_n} \;
    & \sys{\sys{%
	    & \mbf{j_{n+1}} \leq \covD \mbf{C_n}\\
	    & \mbf{j_n} \leq \inD \mbf{C_n} \\
	    & \inD \mbf{C_n} \leq   y_{n+1} \V (x_n \W t_n') \\
	    & \mbf{j_n} \leq x_n \W t_n'\\
	    &  y_{n+1} \leq \kappa(\mbf{j_n}) \\
	    &  \leqjD \mbf{j_{n+1}} \W x_n 
	    \leq \kappa(\mbf{j_n}) \\
	     & x_n \W  s_n'
	    \leq \kappa(\mbf{j_n}) \\
	    & \mbf{j_{n+1}} \W 
	    \covD  \niB ( y_{n+1} \V  (\leqjD \mbf{j_{n+1}} \W x_n ) \V (x_n \W  s_n'))\leq \bot\\
    }
  \implies \false } .
\end{align*}
We can apply now the rules $(\distVW)$ and $(\spW)$ to the third inequality,
and the rule $(\spW)$ on the fourth inequality,  and get:
\begin{align*}
    \forall \overline{x},
    \forall y_{n+1},
    \forall \overline{y},
    \forall \mbf{j_{n+1}},
    \mbf{j_n},
    \forall \mbf{C_n} \;
    & \sys{\sys{%
	    & \mbf{j_{n+1}} \leq \covD \mbf{C_n}\\
	    & \mbf{j_n} \leq \inD \mbf{C_n} \\
	    & \inD \mbf{C_n} \leq   y_{n+1} \V x_n  \\
	    & \inD \mbf{C_n} \leq   y_{n+1} \V  t_n' \\
	    & \mbf{j_n} \leq x_n \\
	    & \mbf{j_n} \leq  t_n'\\
	    &  y_{n+1} \leq \kappa(\mbf{j_n}) \\
	    &  \leqjD \mbf{j_{n+1}} \W x_n 
	    \leq \kappa(\mbf{j_n}) \\
	     & x_n \W  s_n'
	    \leq \kappa(\mbf{j_n}) \\
	    & \mbf{j_{n+1}} \W 
	    \covD  \niB ( y_{n+1} \V  (\leqjD \mbf{j_{n+1}} \W x_n ) \V (x_n \W  s_n'))\leq \bot\\
    }
  \implies \false } .
\end{align*}
By lemma \ref{lem : double ackermann} with the following instantiations
$$t:=t_n' , \quad s:=s_n' , \quad \nomj := \mbf{j_{n+1}} , \quad \nomk:= \mbf{j_n}
 , \quad \mbf{C}:=\mbf{C_n} , \quad x:=x_n,$$
the quasi-inequality above is equivalent to the following quasi-inequality:

\begin{align*}
    \forall x_{n-1}, ..., x_0,
    \forall \overline{y},
    \forall \mbf{j_{n+1}},
    \forall \mbf{C_n} \;
    & \sys{\sys{%
	    & \mbf{j_{n+1}} \leq \covD \mbf{C_n} \\
	    & \mbf{j_{n}} \leq \inD \mbf{C_n} \\
	    & \mbf{j_{n}} \leq t_n'\\
	    & \cl ( \inD \mbf{C_n}  \smallsetminus  \mbf{j_{n}} ) \leq \kappa(\mbf{j_{n}})\\
	    & \leqjD \mbf{j_{n+1}} \W \leqjD \mbf{j_{n}}  \leq \kappa(\mbf{j_{n}}) \\
	    & \leqjD \mbf{j_{n}} \W s_n'\leq \kappa(\mbf{j_{n}}) \\
	    & \mbf{j_{n+1}} \W 
	    \covD  \niB ( \cl ( \inD \mbf{C_n}  \smallsetminus  \mbf{j_{n}} ) 
	    \V (\leqjD \mbf{j_{n+1}} \W \leqjD \mbf{j_{n}} )
	    \V  (\leqjD \mbf{j_{n}} \W s_n'))\leq \bot\\
    }
  \implies \false } .
\end{align*}
where $\cl$ abbreviates the composition $\covD \niB \leqjD$.
By applying the rule $\MinCovB$ bottom to top, the quasi-inequality above can be equivalently rewritten as follows
\begin{align*}
    \forall x_{n-1}, ..., x_0,
    \forall \overline{y},
    \forall \mbf{j_{n+1}},
    \forall \mbf{C_n} \;
    & \sys{\sys{%
	    & \mbf{j_{n+1}} \leq \covD \mbf{C_n} \\
	    & \mbf{j_{n}} \leq \inD \mbf{C_n} \\
	    & \mbf{j_{n}} \leq t_n'\\
	    & \leqjD \mbf{j_{n+1}} \W \leqjD \mbf{j_{n}}  \leq \kappa(\mbf{j_{n}}) \\
	    & \leqjD \mbf{j_{n}} \W s_n'\leq \kappa(\mbf{j_{n}}) \\
	    & \mbf{j_{n+1}} \W 
	    \covD  \niB ( \cl ( \inD \mbf{C_n}  \smallsetminus  \mbf{j_{n}} ) 
	    \V (\leqjD \mbf{j_{n+1}} \W \leqjD \mbf{j_{n}} )
	    \V  (\leqjD \mbf{j_{n}} \W s_n'))\leq \bot\\
    }
  \implies \false } .
\end{align*}
By applying Lemma \ref{lem : MinCovD}, we get:
\begin{align*}
    \forall x_{n-1}, ..., x_0,
    \forall \overline{y},
    \forall \mbf{j_{n+1}},
    \forall \mbf{C_n} \;
    & \sys{\sys{%
	    & \mbf{j_{n+1}} \leq \covD \mbf{C_n} \\
	    & \mbf{j_{n}} \leq \inD \mbf{C_n} \\
	    & \mbf{j_{n}} \leq t_n'\\
	    & \leqjD \mbf{j_{n+1}} \W \leqjD \mbf{j_{n}}  \leq \kappa(\mbf{j_{n}}) \\
	    & \leqjD \mbf{j_{n}} \W s_n'\leq \kappa(\mbf{j_{n}}) \\
    }
  \implies \false } .
\end{align*}
By applying $(\leqAtom)$ and $(\AtCoat)$ to the fourth inequality, 
and $(\leqAtom)$ to the last inequality, we get:
\begin{align*}
    \forall x_{n-1}, ..., x_0,
    \forall \overline{y},
    \forall \mbf{j_{n+1}},
    \forall \mbf{C_n} \;
    & \sys{\sys{%
	    & \mbf{j_{n+1}} \leq \covD \mbf{C_n} \\
	    & \mbf{j_{n}} \leq \inD \mbf{C_n} \\
	    & \leqjD \mbf{j_{n+1}} \W \mbf{j_{n}}\leq \bot \\
	    & \mbf{j_{n}} \leq t_n'\\
	    &  s_n'\leq \kappa(\mbf{j_{n}}) \\
    }
  \implies \false } .
\end{align*}
Notice that the system above consists of a set of pure inequalities 
and a set of inequalities of the exact shape to which the induction hypothesis applies. Since a run of ALBA does not depend on the presence of side pure inequalities, the induction hypothesis implies that ALBA$^l$ succeeds on the system above, and outputs the pure quasi-inequality below, as required:
\begin{align*}
    \forall \mbf{j_{n+1}}, \mbf{j_n},... \mbf{j_0},
    \forall \mbf{C_n}, \mbf{C_{n-1}},...\mbf{C_0}\;
    & \sys{\sys{%
        & \mbf{j_{n+1}} \leq \covD \mbf{C_{n}} \\
	    & \mbf{j_{n}} \leq \inD \mbf{C_{n}}\\
	    & \leqjD\mbf{j_{n+1}} \W \mbf{j_{n}} \leq \bot\\
	    & \mbf{j_n} \leq \covD \mbf{C_{n-1}} \\
	    & \mbf{j_{n-1}} \leq \inD \mbf{C_{n-1}}\\
	    & \leqjD\mbf{j_n} \W \mbf{j_{n-1}} \leq \bot\\
	    & \ldots\\
	    & \mbf{j_1} \leq \covD \mbf{C_{0}} \\
	    & \mbf{j_{0}} \leq \inD \mbf{C_{0}}\\
	    & \leqjD\mbf{j_1} \W \mbf{j_{0}} \leq \bot \\
    }%
  \implies \false } .
\end{align*}

\end{proof}

The lemma below proves the soundness of an Ackermann-type rule for the elimination of non-elementary variables which however cannot be explained in terms of Ackermann principles.

\begin{lemma}
\label{lem : double ackermann}
Let $t$ and $s$ be monotone $\mc{L}^+$-terms such that $x,y\notin Var(t)$. 
For every closed model $\mb{M} = (\mb{E}_L, v)$ such that $\mb{E}_L = ( J(L), \P J(L), \cov,\ni, \leq_J  )$ is the enriched two-sorted frame 
associated with
some finite lattice $L$ 
(cf.\ Definition \ref{def : two sorted frame of a lattice}),
$$ \mb{M} \Vdash (S1) \quad \tiff \quad \mb{M} \Vdash (S2),$$
where
\begin{align*}
    (S1) := \exists x \, \exists y\,
    \exists \nomj \,\exists \nomk \,
     \exists \mbf{C} \;
    \sys{%
	    & \nomj \leq \covD \mbf{C} \\
	    & \nomk \leq \inD \mbf{C} \\
	    & \inD \mbf{C}  \leq y \V x\\
	    & \inD \mbf{C}  \leq y \V t\\
	    & \nomk \leq x\\
	    & \nomk \leq t\\
	    & y \leq \kappa(\nomk)\\
	    & \leqjD \nomj\W x  \leq \kappa(\nomk) \\
	    & x \W s\leq \kappa(\nomk) \\
	    & \nomj \W 
	    \covD  \niB ( y \V  (\leqjD \nomj \W x ) \V (x \W s))\leq \bot\\
    },%
\end{align*}
\begin{align*}
(S2) :=  \exists \nomj \, \exists \nomk \,
     \exists \mbf{C}  \; \sys{%
	    & \nomj \leq \covD \mbf{C} \\
	    & \nomk \leq \inD \mbf{C} \\
	    & \nomk \leq t\\
	    & \cl ( \inD \mbf{C}  \smallsetminus  \nomk ) \leq \kappa(\nomk)\\
	    & \leqjD \nomj\W \leqjD \nomk  \leq \kappa(\nomk) \\
	    & \leqjD \nomk \W s\leq \kappa(\nomk) \\
	    & \nomj \W 
	    \covD  \niB ( \cl ( \inD \mbf{C}  \smallsetminus  \nomk ) 
	    \V (\leqjD \nomj\W \leqjD \nomk )
	    \V (\leqjD \nomk \W s))\leq \bot\\
    },%
\end{align*}
and $\cl(\phi)$ denotes $\covD \niB \leqjD \phi$.
\end{lemma}
\begin{proof}
Assume that  the conjunction of the inequalities in $(S1)$ holds  under $v$. 
Let $v'$ be the $(x,y)$-variant of $v$ such that 
$v'(x) = \leqjD v(\nomk) $ and $v'(y) = \cl(  v(\inD \mbf{C}  \smallsetminus  \nomk ))$.
Since the assignment $v$ is closed,
for any $z \in AtProp \smallsetminus \{x,y\}$, the set $v'(z) =v(z)$ is closed. 
By definition, $v'(y)$ is closed, and $v'(x)$ is closed because for any finite lattice and any $k \in J(L)$, the downset ${\downarrow_{\leq_J}}k$ is a closed set (cf.\ Lemma \ref{lem : closure down j}). 
Thus $v'$ is a closed assignment. 
In addition, $v'(x) \subseteq v(x)$. Indeed, 
the assumption that $\nomk \leq x$ holds under $v$ and $v(x)$ being  closed, hence a downset, imply that ${\downarrow_{\leq_J}} v(\nomk) \subseteq v(x)$, hence we have: 
$$v'(x) = \leqjD v(\nomk) = {\downarrow_{\leq_J}} v(\nomk) \subseteq v(x) .$$

The first, second and third inequalities in $(S2)$ hold under $v'$ since they do not contain the variables $x$ and $y$ and coincide with the first, second and sixth inequalities in $(S1)$,
which by assumption hold under $v$.
The satisfaction of the fifth and sixth inequalities in $(S2)$ under $v'$ is implied by monotonicity, since 
the eighth and ninth  inequalities in $(S1)$ are satisfied under $v$, and since $v'(x) \subseteq v(x)$.
It remains to show that  the fourth and seventh inequalities in $(S2)$ hold under $v'$.
Let $j, k \in J(L)$ and $C \subseteq J(L)$ such that 
$v(\nomj)=\{ j\}$, $v(\nomk)=\{ k\}$ and
$v(\mbf{C}) =  \{ C \}$.
The assumption that 
$\nomj \leq \covD \mbf{C} $ and
$ \nomk \leq \inD \mbf{C}$ hold under $v$ imply that
$C \in \mc{M}(j) $ and $k\in C$.
Hence, $\inD \mbf{C} = C$.
By Lemma \ref{lem : minimal cover closure}.2, 
$k \notin \overline{{\downarrow_{\leq_J}}( C \smallsetminus k)}$. 
Hence,
$$v'(\cl ( \inD \mbf{C}  \smallsetminus  \nomk )) = 
\overline{{\downarrow_{\leq_J}}( C \smallsetminus k)}
\subseteq J(L) \smallsetminus k = v'(\kappa(\nomk)).$$
Thus the fourth inequality in $(S2)$ holds under $v'$.
As to the last inequality, it follows directly from the satisfaction of the previous inequalities under $v'$ and  Lemma \ref{lem : MinCovD}.

\medskip

Let us prove the converse implication. Assume that  the conjunction of the inequalities in $(S2)$ holds  under $v$.  
Let $v'$ be the $(x,y)$-variant of $v$ such that 
$v'(x): = \leqjD v(\nomk )$ and $v'(y) := \cl ( \inD v(\mbf{C})  \smallsetminus  v(\nomk ) )$. 
The first, second and sixth inequalities in $(S1)$ hold under $v'$ since they do not contain the variables $x$ and $y$ and coincide with the first, second and third inequalities in $(S2)$,
which by assumption hold under $v$.
Since $v'(x) = v'(\leqjD \nomk) = {\downarrow_{\leq_J}} v(\nomk)$, the fifth inequality is satisfied under $v'$.
 The satisfaction under $v'$ of the eighth, ninth and tenth inequalities in $(S1)$ immediately follows from the satisfaction of the fifth, sixth and seventh inequalities in $(S2)$ respectively and the definition of $v'$. 

It remains to be shown that the third, fourth and seventh inequalities in $(S1)$ hold under $v'$.
Let $j, k \in J(L)$ and $C \subseteq J(L)$ such that 
$v(\nomj)=\{ j\}$, $v(\nomk)=\{ k\}$ and
$v(\mbf{C}) =  \{ C \}$.
The satisfaction of the first and second inequalities in $(S1)$ under $v'$ imply that $C\in \M(j)$ and $k\in C$, which imply by 
Lemma \ref{lem : minimal cover closure}.2, that $k \notin v'(y)$. This implies that the seventh inequality in $(S1)$ is satisfied under $v'$.
By definition of $v'$ and of the closure, $$v'(\inD \mbf{C}  \smallsetminus  \nomk ) \subseteq \cl ( v'(\inD \mbf{C}  \smallsetminus  \nomk )) = v'(y).$$ In addition, by the satisfaction of the fifth and sixth inequalities in $(S1)$ under $v'$, we have that  $k\in v'(x)$ and $k \in v'(t)$. Hence $$\inD v'(\mbf{C}) = (\inD v'(\mbf{C})  \smallsetminus  k) \cup \{k\} \subseteq v'(y) \cup v'(x)$$ and $$\inD v'(\mbf{C})  = (\inD v'(\mbf{C} ) \smallsetminus  k) \cup \{k\}\subseteq v'(y) \cup v'(t).$$ 
This finishes the proof that the third and fourth inequalities $(S1)$ hold under the closed assignment $v'$. 
\end{proof}

\end{document}